\documentclass[twoside,11pt]{article}

\PassOptionsToPackage{numbers, compress}{natbib}

\usepackage[preprint]{jmlr2e}

\usepackage{verbatim}
\usepackage{algorithm}
\usepackage{algpseudocode} 
\usepackage{color,graphicx}  
\usepackage{amsfonts}   % if you want the fonts
\usepackage{amsmath}
\usepackage{amssymb}
\newtheorem{thm}{Theorem}[section]

\newtheorem{cor}{Corollary}[section]

\newtheorem{ex}{Example}[section]
\newtheorem{rel}{Relation}[section]
\newtheorem{rmk}{Remark}[section]
\newtheorem{ass}{Assumption}[section]
\newcommand{\RR}{\mathbb{R}}      % for Real numbers
\newcommand{\vecc}{\boldsymbol}

\jmlrheading{-}{2021}{-}{-}{2/21}{-}{Soon Hoe Lim}

\ShortHeadings{Anomalous Thermodynamics  in Homogenized Generalized Langevin Systems}{Soon Hoe Lim}
\firstpageno{1}

\begin{document}

\title{Anomalous Thermodynamics \\ in  Homogenized Generalized Langevin Systems}

\author{\name Soon Hoe Lim 
 \email soon.hoe.lim@su.se \\
       \addr Nordita\\
       KTH Royal Institute of Technology and Stockholm University \\
       Stockholm 106 91, Sweden 
       }

\maketitle

\begin{abstract}
We study functionals, such as heat and work, along trajectories of a class of multi-dimensional generalized Langevin systems in various limiting situations that correspond to different level of homogenization. These are the situations where one or more of the inertial time scale(s), the memory time scale(s) and the noise correlation time scale(s) of the systems are taken to zero. We find that, unless one restricts to special situations that do not break  symmetry of the Onsager matrix associated with the fast dynamics, it is generally not possible to express the effective evolution of these functionals solely in terms of trajectory of the homogenized process describing the system dynamics via the widely adopted Stratonovich convention. In fact, an anomalous term is often needed for a complete description, implying that convergence of these functionals needs more information than simply the limit of the dynamical process. We trace the origin of such impossibility  to area anomaly, thereby linking the symmetry breaking and area anomaly. This hold important consequences for many nonequilibrium systems that can be modeled by generalized Langevin equations. Our convergence results hold in a strong pathwise sense.
\\

\end{abstract}

\tableofcontents

\section{Introduction} \label{sect_intro}

We consider a class of non-Markovian Langevin equations, whose coefficients are possibly state-dependent, describing the dynamics of a particle moving in a force field and interacting with the environment.  The evolution of the particle's position,  $\vecc{x}_{t} \in \RR^d$, $t \geq 0$, is given by the solution to the following stochastic integro-differential equation (SIDE) \citep{2019Lim}:
\begin{equation} \label{genle}
 m \ddot{\vecc{x}}_{t} =  \vecc{F}(t, \vecc{x}_{t})  -\vecc{\gamma}_0(\vecc{x}_t) \dot{\vecc{x}}_t  -  \vecc{g}(\vecc{x}_{t}) \int_{0}^{t} \vecc{\kappa}(t-s) \vecc{h}(\vecc{x}_{s}) \dot{\vecc{x}}_{s} ds + \vecc{\sigma}_0(\vecc{x}_t) \vecc{\eta}_t   + \vecc{\sigma}(\vecc{x}_t) \vecc{\xi}_t, 
\end{equation}
with the initial conditions (here the initial time is chosen to be  $t=0$): 
\begin{equation}
\vecc{x}_{0} = \vecc{x}, \ \ \dot{\vecc{x}}_{0} = \vecc{v}.
\end{equation}

In the SIDE \eqref{genle}, overdot denotes derivative with respect to time $t$, $m > 0$ is the mass of the particle, the matrix-valued functions $\vecc{g}: \RR^{d} \to \RR^{d \times q}$, $\vecc{h} : \RR^{d} \to \RR^{q \times d}$, $\vecc{\sigma}: \RR^{d} \to \RR^{d \times r}$, $\vecc{\gamma}_0: \RR^d \to \RR^{d \times d}$ and $\vecc{\sigma}_0: \RR^d \to \RR^{d \times b}$  are the coefficients of the equation, and $\vecc{F} :\RR^+ \times \RR^{d} \to \RR^{d}$ is a force field acting on the particle. Here $d$, $q$, $r$ and $b$ are, possibly distinct, positive integers. Here and throughout the paper, the superscript $^T$ denotes transposition of matrices or vectors and $E$ denotes mathematical expectation. The SIDE \eqref{genle} can be viewed as a Newton's equation of motion (i.e., $m \ddot{\vecc{x}}_t = \vecc{F}(t,\vecc{x}_t)$) with additional forcing terms to be described in the following. 

The second and third term on the right hand side of \eqref{genle} represent the drag experienced by the particle. This drag is modeled by a sum of two deterministic damping terms of different nature. The second term, proportional to the particle's velocity, models instantaneous damping. On the other hand, the third term, involving an integral over the particle's past velocities with the kernel $\vecc{\kappa}(t-s)$, describes non-instantaneous, distributed delayed, damping due to the back-action effects of the environment up to current time.  The matrix-valued function $\vecc{\kappa}: \RR \to \RR^{q \times q}$ is called a memory function and it decays sufficiently fast at infinities. 

The forth and fifth term on the right hand side of \eqref{genle} represent two stochastic forcings (noises) of different nature imparted to the particle. They are $\vecc{\sigma}_0(\vecc{x}_t) \vecc{\eta}_t$, which is a Gaussian white noise, and $\vecc{\sigma}(\vecc{x}_t)\vecc{\xi}_t$,  which is a Gaussian colored noise, both of which are possibly multiplicative. Here the process $\vecc{\eta}_t$ represents a $b$-dimensional white noise, and $\vecc{\xi}_{t}$ is a $r$-dimensional mean zero stationary Gaussian process with the covariance function $\vecc{R}(t) = E[\vecc{\xi}_t \vecc{\xi}_0^T]$. The two noise processes are mutually independent.  The initial conditions $\vecc{x}$ and $\vecc{v}$ are random variables independent of  the noise process $\{(\vecc{\eta}_t, \vecc{\xi}_{t}) : \ t \geq 0\}$. Precise definition and assumptions, as well as physical motivation, for the memory function and the noise processes  will be given in Section \ref{sect_gles}. 

%We will focus on the physical implications of our results concerning the SIDE and therefore we will bypass the rigorous details throughout the paper. In particular, we do not spell out the full rigorous assumptions concerning the SIDE here, unlike the previous work.  Here and throughout the paper, the superscript $^T$ denotes transposition of matrices or vectors and $E$ denotes expectation.

Therefore, \eqref{genle} is a {\it generalized Langevin equation} (GLE), containing the Langevin-Kramers equation studied in \citep{hottovy2015smoluchowski} (by setting $\vecc{h}$ and $\vecc{\sigma}$ to zero) and the GLE studied in \citep{Lim2018} (by setting $\vecc{\gamma}_0$ and $\vecc{\sigma}_0$ to zero) as special cases. The most basic form of GLE, which is a special case of \eqref{genle}, was first introduced by Mori in \citep{mori1965transport} and subsequently used to model many systems in statistical and biological physics \citep{ermak1978brownian}. The GLE has attracted increasing attentions in recent years, due to its successful application in modeling anomalously diffusing systems, active matter systems and many other nonequilibrium systems \citep{goychuk2012viscoelastic,lysy2016model,gottwald_crommelin_franzke_2017,sevilla2018non}.

We remark that GLEs of the form \eqref{genle}, despite being more general in the above sense, are still not the most general ones. Depending on modeling details (for instance, the form of the coupling among various degrees of freedom), one may need to add other forces such as a Basset force (to account for the effect of hydrodynamic backflow \citep{fodor2015generalized}) in the GLEs, or consider GLEs for a set of reaction coordinates/gross variables instead, in which case the resulting GLEs may feature renormalization of bare potential fields, resulting in a potential of mean force (see Section II.B in \citep{hanggi1990reaction} or the recent paper \citep{RevModPhys.92.041002} and the references therein). While it is important to keep in mind of these more general models, we will not study them in this paper.

One particular instance, of important relevance in statistical mechanics, that we will revisit often  is when the coefficients and/or functions defining the GLE \eqref{genle} are related in the following way.

\begin{rel} {\it Fluctuation-dissipation relations (FDRs).} \label{ass_fdr} 
\begin{itemize}
\item[(a)] $\vecc{\sigma}_0 \vecc{\sigma}_0^T = \vecc{\gamma}_0$ (i.e. the fluctuation-dissipation relation of the first kind holds);
\item[(b)] $\vecc{\kappa}(t) = \vecc{R}(t)$ and $\vecc{g} = \vecc{h}^T = \vecc{\sigma}$ (i.e. the fluctuation-dissipation relation of the second kind holds).
\end{itemize}
\end{rel}

It turns out that the GLE \eqref{genle}, with $\vecc{\gamma}_0$ and $\vecc{\sigma}_0$ zero and satisfying (b) in Relation \ref{ass_fdr}, can be derived from a microscopic Hamiltonian model (Kac-Zwanzig or Caldeira-Leggett type) for a small system interacting with a heat bath, or via the Mori-Zwanzig projection approach. See, for instance, Appendix A in \citep{Lim2018} or \citep{hanggi1997generalized,zwanzig1973nonlinear,rey2006open,leimkuhler2018ergodic}. In this case, there will be proportionality constants, containing the temperature of the heat bath as a parameter, in the fluctuation-dissipation relations. Since these constants could be absorbed into  $\vecc{g}$, $\vecc{h}$ or $\vecc{\sigma}$, we choose not to include them explicitly in Assumption \ref{ass_fdr}. Lastly, we remark that the term $-\vecc{\gamma}_0(\vecc{x}_t) \vecc{v}_t$ (when $\vecc{\gamma}_0$ is non-zero) could  be used to model forces of different nature acting on the particle, in particular when $\vecc{\gamma}_0$ is not positive definite (and therefore cannot model a damping term) -- see Example \ref{eg_3}. Throughout this paper, $\vecc{\gamma}_0$ is either zero or non-zero, in which case it is either positive definite or not positive definite. 

%In this paper we will mostly be interested in the case where the instantaneous damping and white noise term are absent, i.e. $\vecc{\gamma}_0 = \vecc{0}$ and $\vecc{\sigma}_0 = \vecc{0}$, in the GLE \eqref{genle}.

%motivation of paper
There are numerous studies focusing on  asymptotic analysis and model reduction of GLEs, aiming to justify the use of low-dimensional phenomenological equations such as the Langevin-Kramers equations and the overdamped Langevin equations for modeling of statistical systems. See, for instance, \citep{ottobre2011asymptotic, Lim2018, nguyen2018small}. There are also many works studying asymptotics of functionals along trajectory of these phenomenological equations \citep{celani2012anomalous, bo2017multiple,ge2018anomalous,PhysRevE.98.052105,Birrell2018}. On the other hand, to our best knowledge works performing asymptotic analysis of functionals along trajectory of generalized Langevin systems, in particular for functionals appearing in stochastic thermodynamics of GLEs, are scarce. 

In this paper we  present a comprehensive multiple time scales analysis (homogenization) of these functionals, as well as of the GLE dynamics, in various limiting situations.  The main goal is {\it to apply the multiscale analysis  to investigate the issue of discretization choice for a class of stochastic integrals appearing in stochastic thermodynamics}. This issue concerns with justification (or not) of the widespread use of Stratonovich convention (midpoint discretization) for defining functionals, such as heat and work, along trajectories of these phenomenological models, used in deriving the  law of energy balance in the energetics literature \citep{sekimoto2010stochastic, seifert2012stochastic,van2013stochastic,peliti2021stochastic}. 

From mathematical viewpoint, the Stratonovich choice of discretization guarantees the vector fields involved transform under a change of coordinates \citep{chetrite2008fluctuation} and is therefore suitable for formulation of coordinate-free SDEs on manifolds. More importantly, from physical point of view, if the work functional is defined as a Stratonovich product, then the exponential of negative entropy production is  the ratio of probabilities of 
time-reversed and forward paths \citep{seifert2005entropy} and is an exponential martingale \citep{chetrite2011two},  leading to a simple differential equation for the entropy production and many interesting mathematical properties \citep{pigolotti2017generic,yang2020unified}. Had the work not been defined with the Stratonovich product, then these nice properties would not hold. However, the Stratonovich choice needs to be carefully justified at a more fundamental level, for instance by taking a GLE as starting point for analysis, in which case the functionals (stochastic integrals) along the phase-space trajectories are uniquely defined (i.e. their discretization is free of ambiguities). Performing homogenization on these functionals allows us to find out its limiting expression in the considered limit. This limiting expression is then compared to the functional defined along the trajectory of the limiting dynamics.

In our previous contribution in \citep{Bo19}, we have shown that for systems in which noise correlation is shorter-lived than inertia (usually the case for microscopic colloids in water at room temperature) the correct discretization for these functionals is Stratonovich -- this is the result obtained by performing a Markovian limit first and then the small mass limit. This result holds under the conditions that (i) the processes which generate the colored noise are equilibrium ones, and (ii) in the small mass limit the velocity degrees of freedom reach an equilibrium distribution with the local temperature (this holds when the fluctuation-dissipation relation is obeyed). For systems that violate these conditions, the interpretation of the (limiting) functionals is less immediately clear.  The main motivation and contribution of this paper is, in fact, to investigate and identify the limiting behavior of these functionals beyond the aforementioned setting via a systematic multiscale analysis considering different hierarchies of the time scales involved. The results obtained in this paper not only recover our earlier results in \citep{Bo19}, but also give new results and uncover interesting insights in more general settings. We emphasize that the present paper has a rather applied flavor, i.e., it focuses on applying  the general homogenization theorems obtained in our earlier works \citep{2019Lim,Lim2018}  to shed light on the above issue in the field of stochastic thermodynamics, rather than presenting novel mathematical techniques and proofs for homogenization. Moreover, the notion of stochastic L\'evy areas is, for the first time, connected to thermodynamic quantities such as work done on physical systems.

%It can also be seen as a review paper, collecting present and recent   homogenization results for generalized Langevin systems in one place.

%More importantly, to cater to a wider audience, we organize the statements of our results into two different levels of display at the cost of lengthiness and redundancy: results in the first level bypasses all the technical assumptions and details, focusing on mostly the physical situations and the central messages to be conveyed, whereas results in the second level are mathematically rigorous convergence statements, obtained via a set of technical assumptions, covering  more general situations than the ones considered in Section \ref{sect_limit}.

%strong pathwise convergence

This paper is organized as follows. In Section \ref{sect_gles}, we define the class of GLE models to be studied in this paper. We give three examples, of relevance in applications to study nonequilibrium systems, of these models in Appendix \ref{sect_eg}. In Section \ref{sect_funct}, we motivate and introduce a class of  functionals along trajectories of the GLE. 
In Section \ref{sect_gen_homog}, we study homogenization for a class of SDE systems with state-dependent coefficients and their functionals. The convergence results are obtained in a strong pathwise sense. They follow from applications of the homogenization theorem  in Appendix \ref{sect_generalhomogthm}. We discuss the mathematical implications of these results, in particular we link symmetry breaking of the Onsager matrix associated to the fast dynamics and area anomaly.  In Section \ref{app_levy}, we illustrate and discuss this link in the context of a Brownian particle in a magnetic field to build some intuition on area anomaly before moving on to study the more general situations of GLEs. Section \ref{sect_limit} contains the main contributions of the paper. There, building on the results in Section \ref{sect_gen_homog}, we study homogenization for generalized Langevin dynamics as well as the functionals introduced in Section \ref{sect_funct}.  We then discuss the conditions under which a Stratonovich functional is recovered for various limiting situations, as well as the consequences due to interplay between the symmetry breaking  and area anomaly.  We conclude the paper in Section \ref{sect_concl}. 

% All other supplementary materials and technical details are deferred to the appendices. 

%In Section \ref{sect_appl}, we apply the general results to study  GLE model for a number of nonequilibrium systems. 

\section{Generalized Langevin Equations (GLEs)} \label{sect_gles}

In this section we define our GLE models, following closely the notation in \citep{Lim2018}. In the GLE \eqref{genle}, the memory function $\vecc{\kappa}: \RR \to \RR^{q \times q}$ is taken to be {\it Bohl}, i.e. the matrix elements of $\vecc{\kappa}(t)$ are finite linear combinations of the functions of the form $t^k e^{\alpha t} \cos(\omega t)$ and $t^k e^{\alpha t} \sin(\omega t)$, where $k$ is an integer and $\alpha$ and $\omega$ are real numbers. For properties of Bohl functions, we refer to Chapter 2 of \citep{trentelman2002control}.  The noise process $\vecc{\xi}_{t}$ is a $r$-dimensional mean zero stationary real-valued Gaussian vector process having a Bohl covariance function, $\vecc{R}(t):=E \vecc{\xi}_t \vecc{\xi}_0^T = \vecc{R}^T(-t) $, and, therefore, its spectral density, $\vecc{S}(\omega) := \int_{-\infty}^{\infty} \vecc{R}(t) e^{-i\omega t} dt$, is a rational function \citep{willems1980stochastic}.

Note that the Gaussian process $\vecc{\xi}_t$ which drives the SIDE \eqref{genle} is not assumed to be Markov.  The assumptions we made on its covariance will allow us to present it as a projection of a Markov process in a (typically higher-dimensional) space.  This approach, which originated in stochastic control theory \citep{kalman1960new}, is called {\it stochastic realization}. We describe $\vecc{\kappa}(t)$ and $\vecc{\xi}_t$ in detail below.

Let $\vecc{\Gamma}_1 \in \RR^{d_1 \times d_1}$, $\vecc{M}_1 \in \RR^{d_1 \times d_1}$, $\vecc{C}_1 \in \RR^{q \times d_1}$, $\vecc{\Sigma}_1 \in \RR^{d_1 \times q_1}$, $\vecc{\Gamma}_2 \in \RR^{d_2 \times d_2}$, $\vecc{M}_2 \in \RR^{d_2 \times d_2}$, $\vecc{C}_2 \in \RR^{r \times d_2}$, $\vecc{\Sigma}_2 \in \RR^{d_2 \times q_2}$ be constant matrices, where $d_1,d_2,q_1,q_2$, $q$ and $r$ are positive integers. 
In this paper, we study the class of SIDE \eqref{genle}, with the memory function defined in terms of the triple $(\vecc{\Gamma}_1,\vecc{M}_1,\vecc{C}_1)$ of matrices as follows:
\begin{equation} \label{memory_realized}
\vecc{\kappa}(t)=\vecc{C}_1e^{-\vecc{\Gamma_1}|t|}\vecc{M}_1\vecc{C}_1^T.
\end{equation}
The covariance of the stationary Gaussian noise process $\vecc{\xi}_t$ will be expressed in terms of the triple $(\vecc{\Gamma}_2,\vecc{M}_2,\vecc{C}_2)$.  More precisely, we define it as:
\begin{equation} \label{noise}
\vecc{\xi}_t = \vecc{C}_2 \vecc{\beta}_t,\end{equation}
where $\vecc{\beta}_t$ is the solution to the It\^o SDE:
\begin{equation} \label{realize}
d\vecc{\beta}_t = -\vecc{\Gamma}_2\vecc{\beta}_t dt + \vecc{\Sigma}_2 d\vecc{W}^{(q_2)}_t,
\end{equation}
with the initial condition, $\vecc{\beta}_0$, normally distributed with zero mean and covariance  $\vecc{M}_2$. Here, $\vecc{W}^{(q_2)}_t$ denotes a $q_2$-dimensional Wiener process and is independent of $\vecc{\beta}_0$. 
%Throughout the paper the dimension of the Wiener process will be specified by the superscript in this way. 

For $i=1,2$, the matrix $\vecc{\Gamma}_i$ is {\it positive stable}, i.e. all its eigenvalues have positive real parts and $\vecc{M}_i = \vecc{M}_i^T > 0$ satisfies the following Lyapunov equation:
\begin{equation} 
\vecc{\Gamma}_i \vecc{M}_i+\vecc{M}_i \vecc{\Gamma}_i^T=\vecc{\Sigma}_i \vecc{\Sigma}_i^T.
\end{equation}
It follows from positive stability of $\vecc{\Gamma}_i$ that this equation indeed has a unique solution \citep{bellman1997introduction}. 
The covariance matrix, $\vecc{R}(t) \in \RR^{r \times r}$, of the noise process is therefore expressed in terms of  the matrices $(\vecc{\Gamma}_2,\vecc{M}_2,\vecc{C}_2)$ as follows:
\begin{equation} \label{cov}
\vecc{R}(t)=\vecc{C}_2e^{-\vecc{\Gamma_2}|t|}\vecc{M}_2\vecc{C}_2^T, 
\end{equation}
and so the triple $(\vecc{\Gamma}_2,\vecc{M}_2,\vecc{C}_2)$ completely specifies the probability distribution of $\vecc{\xi}_t$. For concrete examples of noise process that can be realized using the above formalism, see \citep{Lim2018}.

%In this paper, we take as the starting point the realization of the memory function and noise process defined by the triples $(\vecc{\Gamma}_i, \vecc{M}_i, \vecc{C}_i)$ as above, and will not be concerned with the details and algorithms associated with problems such as the spectral factorization problem in the stochastic realization theory. Our consideration, in fact, covers the largest class of systems (i.e. with a Bohl memory function and driving noise process having a rational spectral density) that one can possibly realize in a finite dimensional state space (see the propositions and theorems on page 303-308 of \citep{willems1980stochastic}). One could in principle extend the studies to the case of realizations on infinite dimensional state space, in which case the matrices become operators on a Hilbert space and a richer class \footnote{Indeed, every wide sense separable Gaussian process admits a stochastic realization on a suitable Hilbert space (see Proposition 4 in \citep{lindquist1978minimal}).} of memory functions and noise processes can be realized, but this will not be investigated here. 

Physically, the choice of the matrices  $\vecc{\Gamma}_2,\vecc{M}_2,\vecc{C}_2$ specifies the characteristic time scales (eigenvalues of $\vecc{\Gamma}_2^{-1}$) present in the environment, introduces the initial state of a stationary Markovian Gaussian noise and selects the parts of the prepared Markovian noise that are (partially) observed, respectively.  In other words, we have assumed that the  noise in the SIDE \eqref{genle} is realized or ``experimentally prepared" by the above triple of matrices \citep{Lim2018}. The triples that specify the memory function in \eqref{memory_realized} and the noise process in \eqref{noise} are unique up to the following transformations:  
\begin{equation} \label{transf_realize}
(\vecc{\Gamma}'_i=\vecc{T}_i \vecc{\Gamma}_i \vecc{T}^{-1}_i, \vecc{M}_i' = \vecc{T}_i \vecc{M}_i \vecc{T}_i^{T}, \vecc{C}'_i =  \vecc{C}_i \vecc{T}_i^{-1}),
\end{equation}
where $i=1,2$ and the $\vecc{T}_i$ are invertible matrices of appropriate dimensions.

%The triple $(\vecc{\Gamma}_2,\vecc{M}_2,\vecc{C}_2)$ above is called a {\it (weak) stochastic realization} of the covariance matrix $\vecc{R}(t)$ in the well established theory of stochastic realization, which is concerned with solving the inverse problem of stationary covariance generation (see \citep{lindquist1985realization,lindquist2015linear}).  Any zero mean stationary Gaussian process, $\vecc{\xi}'_t$, having a Bohl covariance function, can be realized as a projection of a Gaussian Markov process in the above way.

With the above definitions of memory kernel and noise process, the SIDE \eqref{genle} becomes:
\begin{equation} \label{genle_general}
 m \ddot{\vecc{x}}_{t} =  \vecc{F}(t,\vecc{x}_{t}) - \vecc{\gamma}_0(\vecc{x}_t) \dot{\vecc{x}}_t -  \vecc{g}(\vecc{x}_{t}) \int_{0}^{t} \vecc{C}_1 e^{-\vecc{\Gamma}_1(t-s)} \vecc{M}_1\vecc{C}_1^T \vecc{h}(\vecc{x}_{s}) \dot{\vecc{x}}_{s} ds + \vecc{\sigma}_0(\vecc{x}_t) \vecc{\eta}_t +   \vecc{\sigma}(\vecc{x}_{t}) \vecc{C}_2 \vecc{\beta}_{t}, 
\end{equation}
where $\vecc{\beta}_t$ is the solution to the SDE \eqref{realize}. Introducing the auxiliary variable 
\begin{equation}
\vecc{y}_t = \int_0^t e^{-\vecc{\Gamma}_1(t-s)} \vecc{M}_1 \vecc{C}_1^T \vecc{h}(\vecc{x}_s) \vecc{v}_s ds,
\end{equation}
and setting $\vecc{\eta}_t dt=d \vecc{B}_t$, where $\vecc{B}_t \in \RR^b$ is a Wiener process independent of $\vecc{W}_t^{(q_2)}$, 
the SIDE can be cast as the following It\^o SDE system for the Markov process $\vecc{z}_t = (\vecc{x}_t, \vecc{v}_t, \vecc{y}_t, \vecc{\beta}_t) \in \RR^{d \times d \times d_1 \times d_2}$: 
\begin{align}
d\vecc{x}_t  &= \vecc{v}_t dt, \label{gle1} \\ 
m d\vecc{v}_t &= \vecc{F}(t,\vecc{x}_t) dt - \vecc{\gamma}_0(\vecc{x}_t) \vecc{v}_t dt -   \vecc{g}(\vecc{x}_t) \vecc{C}_1 \vecc{y}_t dt + \vecc{\sigma}_0(\vecc{x}_t) d\vecc{B}_t + \vecc{\sigma}(\vecc{x}_t) \vecc{C}_2 \vecc{\beta}_t dt, \\
d\vecc{y}_t &= -\vecc{\Gamma}_1 \vecc{y}_t dt + \vecc{M}_1 \vecc{C}_1^T \vecc{h}(\vecc{x}_t) \vecc{v}_t dt, \\ 
d\vecc{\beta}_t &= -\vecc{\Gamma}_2 \vecc{\beta}_t dt + \vecc{\Sigma}_2 d\vecc{W}_t^{(q_2)}.  \label{gle4}
\end{align}

We refer to Appendix \ref{sect_eg} for three examples of GLE system arising in nonequilibrium statistical mechanics. Several remarks concerning the system \eqref{gle1}-\eqref{gle4} are now in order.

\begin{rmk}
On one hand, $\vecc{z}_t$ is the solution to a hypoelliptic  SDE system of the form 
\begin{equation} \label{hypo_gen}
d\vecc{z}_t = \vecc{a}(t,\vecc{z}_t)dt + \vecc{B}(t,\vecc{z}_t) d\vecc{U}_t,
\end{equation}
where $\vecc{U}_t$ is a Wiener process and $\vecc{B}$ is a matrix-valued function that is not full rank, since the noise does not act in all directions of $\vecc{z}$. Therefore, from mathematical point of view our study of the GLE and functionals along its trajectory can be viewed as study of the above hypoelliptic SDE system \citep{pavliotis2014stochastic} and the associated functionals.
On the other hand, the process $\vecc{r}_t = (\vecc{x}_t, \vecc{v}_t, \vecc{y}_t)$ gives the coordinates of the generalized Langevin system. It is a non-Markov process satisfying an It\^o SDE of the form:
\begin{equation}
d\vecc{r}_t = \vecc{b}(t, \vecc{r}_t) dt + \vecc{\Phi}(\vecc{r}_t) d\vecc{B}_t + \vecc{\Phi}_a(\vecc{r}_t) \vecc{\beta}_t dt,  \label{r_eq}
\end{equation}
where the driving noise consists of a white noise and a Gaussian colored noise. Note that the augmented process $\vecc{z}_t = (\vecc{r}_t, \vecc{\beta}_t)$ is the Markov process solving the SDE \eqref{hypo_gen}. 
\end{rmk}

%We now mention a few special cases of interest. When $\vecc{g} = \vecc{h}^T = \vecc{\sigma}$, i.e. when Assumption \ref{ass_fdr} is satisfied, we have a lower dimensional  SDE, satisfied by the Markov process $(\vecc{x}_t,\vecc{v}_t, \vecc{u}_t)$, where $\vecc{u}_t = \vecc{y}_t - \vecc{\beta}_t$:
%\begin{align}
%d\vecc{x}_t  &= \vecc{v}_t dt, \label{gle11} \\ 
%m d\vecc{v}_t &= \vecc{F}(t,\vecc{x}_t) dt - (\vecc{\sigma}_0 \vecc{\sigma}_0^T)(\vecc{x}_t) \vecc{v}_t dt -   \vecc{\sigma}(\vecc{x}_t) \vecc{C}_1 \vecc{u}_t dt + \vecc{\sigma}_0(\vecc{x}_t) d\vecc{B}_t, \label{gle22} \\
%d\vecc{u}_t &= -\vecc{\Gamma}_1 \vecc{u}_t dt + \vecc{M}_1 \vecc{C}_1^T \vecc{\sigma}^T(\vecc{x}_t) \vecc{v}_t dt - \vecc{\Sigma}_1 d\vecc{W}_t^{(q_2)}. \label{gle33}
%\end{align} 

\begin{rmk}
One could have absorbed the constant matrices $\vecc{C}_i$ into the coefficients $\vecc{\sigma}$, $\vecc{g}$, $\vecc{h}$ but we choose to keep them as parameters for our memory function and colored noise models. The one-dimensional case $(d=1)$ where $\vecc{C}_i = 1$, $\vecc{\Gamma}_i = \alpha_i > 0$, $\vecc{\Sigma}_i = \alpha_i$, $\vecc{M}_i = \alpha_i/2$, for $i=1,2$ (we will drop the boldface when denoting the processes and coefficients in the one-dimensional case -- for instance, $\vecc{x}_t = x_t$, $\vecc{g} = g$, $\vecc{W}_t = W_t$, etc.), follows as a special case. In this case, the memory function and covariance function of the colored noise process are exponentials, with possibly different decay rates $\alpha_i$. 
\end{rmk}

%We will also keep in mind the two distinct cases where the instantaneous damping and white noise term (Markovian force) are either present or absent in the GLE. 

\begin{rmk}
In order to be able to study the GLE as a finite-dimensional Markovian system it is crucial that the memory function and covariance function of the colored noise process be Bohl. In the case where, for instance, these functions decay as a power law, the resulting GLE cannot be studied as a finite-dimensional SDE system and one needs to work in the infinite-dimensional setting \citep{Kupferman2004, glatt2018generalized}. In other words, one needs to introduce an infinite number of auxiliary variables for describing the memory process $\vecc{y}_t$ and noise process $\vecc{\xi}_t$ to make the extended system Markovian. However, our formalism allows us to approximate an arbitrary memory function, such as the ones decaying as a power law (long-range memory), on a finite time scale \citep{siegle2011markovian}. Therefore, our finite-dimensional consideration allows us to cover a sufficiently large class of systems with memory.
\end{rmk}

\section{Functionals Along Trajectories of GLEs} \label{sect_funct}

We are interested in the asymptotic behavior of a class of functionals along the trajectory $(\vecc{r}_t)_{t \geq 0}$, where $\vecc{r}_t = (\vecc{x}_t, \vecc{v}_t, \vecc{y}_t)$\footnote{Since $\vecc{y}_t$ is a functional of $(\vecc{x}_s, \vecc{v}_s)_{0 \leq s \leq t}$, it suffices to consider the trajectory $(\vecc{x}_t, \vecc{v}_t)_{t \geq 0}$ instead of  $(\vecc{x}_t, \vecc{v}_t, \vecc{y}_t)_{t \geq 0}$ }, of the generalized Langevin systems described by \eqref{genle_general} in various limiting situations. These situations are when wide separation of time scales exists in the systems and thereby allowing simplification of the dynamics via elimination of the fast degrees of freedom and description of the system solely in terms of the slow degrees of freedom.  These functionals take the form of:
\begin{equation} \label{gen_func}
\mathcal{F}_t = \int_0^t r(s, \vecc{r}_s) ds + \int_0^t  \vecc{p}(s, \vecc{r}_s) \circ^{?}  d\vecc{r}_s 
%+ \int_0^t \vecc{p}_2(s, \vecc{x}_s, \vecc{v}_s) \cdot d\vecc{v}_s,
\end{equation}
which, in differential form, is:
\begin{equation}
d\mathcal{F}_t = r(t, \vecc{r}_t) dt +  \vecc{p}(t, \vecc{r}_t)  \circ^{?} d\vecc{r}_t, 
\end{equation}
where $\circ^{?}$ denotes the (to be specified) discretization rule defining the stochastic integral in \eqref{gen_func}. Since different discretization rules lead to different properties of the functional, the discretization rule should be assigned in such a way that the physical behavior of the modeled system is captured correctly \citep{hottovy2012noise,farago2014langevin,sokolov2010ito,yang2017hydrodynamic}.  Here and throughout the paper, we are using calligraphic font for denoting a functional.  We emphasize that, in contrast to the case of Langevin-Kramers model, the process $\vecc{r}_t$, being a component of the Markov process $(\vecc{r}_t, \vecc{\beta}_t)$, is generally non-Markov. 

%The above functionals  belong to the more general class of functionals  that are of the form:
%\begin{equation}
%\mathcal{F}_t = \int_0^t r(s, \vecc{z}_s) ds + \int_0^t \vecc{p}(s, \vecc{z}_s) \circ d\vecc{z}_s,
%\end{equation}
%where $\vecc{z}_t$ is the solution to the hypoelliptic SDE \eqref{hypo_gen} and will collect the slow and fast variables in the case of homogenization later.  

We are going to introduce and define a special subclass of functionals \eqref{gen_func} along the trajectory of the GLE \eqref{genle_general} (or equivalently the SDE system \eqref{gle1}-\eqref{gle4}) in the following. These functionals are various thermodynamic functionals of interest arising in stochastic thermodynamics  \citep{seifert2012stochastic} of the GLE. To begin with, we split the force field as  $\vecc{F}(t, \vecc{x}) = -\vecc{\nabla}_{\vecc{x}} U(t,\vecc{x})+\vecc{f}_{nc}(t, \vecc{x})$, where the scalar-valued function $U$ represents a potential and $\vecc{f}_{nc}$ represents a non-conservative external force, driving the system out of equilibrium.  
%Due to the hypoellipticity nature of the SDE system \eqref{hypo_gen}, we expect to encounter some difficulties when trying to define these functionals in a meaningful way. 

When considering these functionals, there are two cases of interest. The first case is the case when $\vecc{\sigma}_0 = \vecc{0}$,  in which case there is no ambiguity in defining the stochastic integral in \eqref{gen_func}.  The second case is when $\vecc{\sigma}_0$ is  non-zero, in which case we need to specify the convention $\circ^{?}$ for the stochastic integral, usually taken to be Stratonovich.  We will consider only the first case here. Therefore, {\it we set $\vecc{\sigma}_0$ to zero from now on unless specified otherwise}, and replace $\circ^{?}$ by $\cdot$ to denote dot product. More precisely, when $\vecc{\sigma}_0$ vanishes (and therefore the corresponding $\vecc{\Phi}$ in \eqref{r_eq} vanishes), the equation for $\vecc{r}_t$ does not contain a white noise term. In this case, the process $\vecc{r}_t$ is more regular than the one in the case of non-vanishing $\vecc{\sigma}_0$ and  the stochastic integral defining $\vecc{r}_t$ is uniquely defined, in particular its properties are independent of the discretization choice. 

%\footnote{We remark that when $\vecc{\sigma}_0$ is non-zero, one can, at least in special cases, derive a fluctuation relation for the GLE. }

We define a {\it heat-like} and {\it work-like} functional along the stochastic trajectory $(\vecc{r}_t)_{t \geq 0}$ as the functional satisfying the following (controlled) differential equations:
\begin{align}
d\mathcal{Q}_t &= \left(- \vecc{g}(\vecc{x}_t)\int_0^t \vecc{\kappa}(t-s) \vecc{h}(\vecc{x}_s) \vecc{v}_s ds + \vecc{\sigma}(\vecc{x}_t) \vecc{\xi}_t - \vecc{\gamma}_0(\vecc{x}_t) \vecc{v}_t \right) \cdot d\vecc{x}_t, \\ 
&=  \int_{0}^t \bigg( m \vecc{v}_s \cdot d\vecc{v}_s - \vecc{F}(s,\vecc{x}_s) \cdot d\vecc{x}_s  \bigg),  \label{heat_def} \\
d\mathcal{W}_t &= \frac{\partial U}{\partial t} dt + \vecc{f}_{nc}(t, \vecc{x}_t)\cdot d\vecc{x}_t \label{work_def}
\end{align}
respectively. The above functionals are free of ambiguities in the discretization procedure and are thus uniquely defined.

We emphasize that, as we discussed in \citep{Bo19}, the functionals above are not, generally and strictly speaking, defining  physical heat and work for the generalized Langevin systems. This emphasis leads to our usage of the terminology ``heat-like'' and ``work-like'' functional instead of heat and work throughout the paper. These heat-like and work-like functionals are rather defined in a manner that ensures  a first law for energy balance  is satisfied as follows. Let us define\footnote{Note that a fluctuating internal energy is by no means  unique, but can assume many different forms which all would give the same ``mean value'', but different higher moments \citep{hanggicom}. As a consequence, \eqref{def_int} is nothing more than a definition.} the internal energy of the system as:
\begin{equation} 
\mathcal{E}_t =\frac{1}{2} m |\vecc{v}_t|^2+U(t,\vecc{x}_t). \label{def_int}
\end{equation}
Then, the above definitions for heat-like and work-like functional are consistent with the first law of stochastic thermodynamics in the sense that the energy $\mathcal{E}_t$ is conserved along individual trajectories. Indeed, using $d\mathcal{E} = m \vecc{v} \cdot d\vecc{v} + dU$, one obtains the law:
\begin{equation}
d\mathcal{E} = d\mathcal{W}+ d\mathcal{Q},
\end{equation}
where $\mathcal{W}$ and $\mathcal{Q}$ are defined in \eqref{heat_def} and \eqref{work_def} respectively, and we use the convention that $\mathcal{Q}<0$ if the heat is transferred or dissipated from the system into the environment. 

%The entropy of the particle is the state function  given by the Shannon entropy: 
%\begin{equation}  S_t^{par}= - k_B \log \rho_S(t, \vecc{x}_t,\vecc{v}_t),
%\end{equation}
%where $\rho_S$ is the marginal probability distribution of the solution to the Fokker-Planck equation associated with the SDE system \eqref{gle1}-\eqref{gle4}. 

Next, we specialize the above definition to the setting where the heat-like and work-like functional become physical heat and work. This is the case where $\vecc{\gamma}_0 = \vecc{0}$, the fluctuation-dissipation relation of the second kind holds, and the colored noise models a heat bath which is in equilibrium at temperature $T$. In this case, the resulting GLE can be derived from a microscopic Hamiltonian model (see an earlier remark in Section \ref{sect_gles}) for a Brownian particle (weakly) interacting with an equilibrium  heat bath at temperature $T$. The {\it thermodynamic entropy produced in the environment}, from an initial state $(\vecc{x}_0,\vecc{v}_0)$ at the initial time to a final state $(\vecc{x}_t,\vecc{v}_t)$ at time $t$, is  defined as:
\begin{align}
\mathcal{S}_t = -\beta \mathcal{Q}_t = \beta \int_{0}^t \bigg(  \vecc{F}(s,\vecc{x}_s) \cdot d\vecc{x}_s - m \vecc{v}_s \cdot d\vecc{v}_s \bigg).  \label{entropy}
\end{align}
where $\beta = 1/k_B T$. It is a measure of irreversibility of the generalized Langevin dynamics. The heat can be interpreted as the change of bath energy over the system trajectory and it is a functional of the system history alone \citep{aurell2018unified}.  In the more general case beyond the above setting, the above definition does not generally define a thermodynamic entropy, and so we are going to simply refer to it as an entropy-like functional. Finally, we emphasize that the integrals defining the dynamical process $\vecc{r}_t$ and functionals $\mathcal{Q}_t$, $\mathcal{R}_t$ here are {\it uniquely defined and will be taken to be the starting point for multiple time scale analysis (homogenization)}, for which (the interpretation of) their limiting expression will be of interest.

\section{Homogenization of Slow-Fast SDE Systems and Their Functionals}
\label{sect_gen_homog}

Asymptotic analysis of functionals along trajectories of  approximating stochastic processes has long histories and is an important tool for stochastic modeling of noisy systems. An important early example comes from the classic work of Wong and Zakai \citep{wong1965convergence}, who considered the limiting behavior of the family of real-valued stochastic integrals $y_n(t) = \int_0^t u(B_n(s)) dB_n(s)$, where $u$ is some sufficiently nice function and $B_n(t)$ is a sequence of sufficiently smooth functions approximating a Wiener process. They found that $y_n(t)$ converges to the Stratonovich integral, $y(t) = \int_0^t u(B(s)) \circ dB(s)$, where $\circ$ denotes Stratonovich product and $B(t)$ is a Wiener process, in the limit as $n \to \infty$. The result holds in one dimension and may fail in higher dimensions, in which case one has additional (anomalous) drift terms due to L\'evy  area correction  \citep{levy1951wiener, ikeda2014stochastic,sussmann1991limits} (see Section 11.7.7 in \citep{pavliotis2008multiscale} for an explicit example).  

Each $y_n(t)$ is a functional along trajectories of the approximating functions $B_n(t)$. In the special case where the fast process $B_n(t)$ satisfies an It\^o SDE, driven by a white noise, the key technique is to embed the functional into a higher dimensional Markov process. The goal is then to determine the limiting behavior of the slow process $y_n(t)$, as components of the Markov process, as $n \to \infty$. In the context of the above example, one has $dz_n(t) = dB_n(t)$, $dy_n(t) = u(z_n(t)) dz_n(t)$, and $B_n(t)$ is a process embedded in a SDE system. If, for instance,  $B_n(t)$ is an integrated Ornstein-Uhlenbeck process, then we have $dB_n(t) = C_n(t) dt$, $dC_n(t) = -\lambda_n C_n(t) dt + \sigma_n dW_t$, where $W_t$ is a Wiener process and $\lambda_n$, $\sigma_n$ are some suitable increasing sequences in $n$. 

We are going to study a  generalization of the above example problem to a class of multi-dimensional diffusion processes. Our setting is sufficiently general to cover all the asymptotic  problems (homogenization) for GLEs and their functionals in this paper. In this section we focus on homogenization in the general setting. Examples and applications in the context of stochastic energetics will be studied and discussed in detail in Section  \ref{app_levy} and Section \ref{sect_limit}.
%We will study the problem using two different approaches: (1) formal perturbative expansion of the infinitesimal generator, and (2) homogenization at the level of sample path. The results obtained via the first approach could be justified using theorems of Kurtz. Although the first approach gives very weak convergence result, it is more widely used and intuitive. The second approach leads to strong convergence result but the analysis is often more difficult to be extended to a more general problem. 

Consider the following family of It\^o SDE systems  for $\vecc{Z}^\epsilon_t = (\vecc{X}^\epsilon_t, \vecc{Y}^\epsilon_t, \vecc{\mathcal{A}}^\epsilon_t, \mathcal{B}^\epsilon_t)\in \RR^n \times \RR^m \times \RR^l \times \RR$:  
\begin{align}
d\vecc{X}^\epsilon_t &=  \vecc{U}_1(t,\vecc{X}^\epsilon_t) \vecc{Y}^\epsilon_t dt + \vecc{u}_1(t,\vecc{X}^\epsilon_t) dt + \tilde{\vecc{\sigma}}(t,\vecc{X}^\epsilon_t) d\tilde{\vecc{W}}_t, \label{eqx} \\ 
\epsilon d\vecc{Y}^\epsilon_t &= - \vecc{U}_2(t,\vecc{X}^\epsilon_t) \vecc{Y}^\epsilon_t dt + \vecc{u}_2(t,\vecc{X}^\epsilon_t) dt + \vecc{\sigma}(t,\vecc{X}^\epsilon_t) d\vecc{W}_t, \label{eqy} \\ 
d \vecc{\mathcal{A}}^\epsilon_t &= \vecc{r}(t, \vecc{X}^\epsilon_t) dt  +  \vecc{P}(t,\vecc{X}^\epsilon_t)  d\vecc{X}^\epsilon_t, \label{eqA}  \\
d\mathcal{B}^\epsilon_t &=  \epsilon \vecc{Y}^\epsilon_t \cdot d\vecc{Y}^\epsilon_t, \label{egB}   
\end{align}
where $\vecc{U}_1: \RR^+ \times \RR^n \to \RR^{n \times m}$,  $\vecc{U}_2: \RR^+ \times \RR^n \to \RR^{m \times m}$,   $\vecc{u}_1: \RR^+ \times \RR^n \to \RR^n$,  $\vecc{u}_2: \RR^+ \times \RR^n \to \RR^m$, $\tilde{\vecc{\sigma}}: \RR^+ \times \RR^n \to \RR^{n \times d_s}$,   $\vecc{\sigma}: \RR^+ \times \RR^n \to \RR^{m \times d_f}$,    $\tilde{\vecc{W}}_t \in \RR^{d_s}$ and $\vecc{W}_t \in \RR^{d_f}$ are independent Wiener processes on a filtered probability space $(\Omega, \mathcal{F}, \mathcal{F}_t, \mathbb{P})$ such that the usual conditions \citep{karatzas2014brownian} hold, $\vecc{r}: \RR^+ \times \RR^n \to \RR^l$, $\vecc{P}: \RR^+ \times \RR^n \to \RR^{l \times n}$,  $\epsilon > 0$ is a small parameter, and $\cdot$ denotes dot product. The variables $\vecc{Z}^\epsilon_t$ model physical processes or states of a system with dimensionless variables. Let $\mathbb{E}$ denote expectation with respect to $\mathbb{P}$.

%Here we have chosen to work with a specific fast process to allow explicit expression for the effective coefficients to be derived later as well as to obtain a strong pathwise convergence result, rather than more general ergodic processes  \citep{pavliotis2008multiscale}.

We take $\mathcal{B}_0^\epsilon = \epsilon |\vecc{Y}_0^\epsilon|^2/2$, so that 
\begin{align}
\mathcal{B}_t^\epsilon &= \mathcal{B}_0^\epsilon + \epsilon \int_0^t \vecc{Y}_s^\epsilon \cdot d\vecc{Y}_s^\epsilon = \mathcal{B}_0^\epsilon + \frac{\epsilon}{2} \int_0^t d \left(|\vecc{Y}_s^\epsilon|^2 \right)  = \frac{\epsilon}{2} |\vecc{Y}_t^\epsilon|^2. \end{align}

%We will also keep in mind the special case when $q = 0$ and $\vecc{p} = \vecc{p}(\vecc{x})$. 

The above systems are variants of the one considered in \citep{bo2013white} (see also \citep{bo2017multiple,bo2014entropy}). All the equations contain fast dynamics but the dynamics in $\vecc{Y}^\epsilon$ is one order of magnitude faster than in $\vecc{X}^\epsilon$, $\vecc{\mathcal{A}}^\epsilon$ and $\mathcal{B}^\epsilon$. Our goal is to eliminate the variable $\vecc{Y}^\epsilon$  in \eqref{eqx}-\eqref{egB} and derive an effective description for the slow process $\vecc{Q}^\epsilon_t = (\vecc{X}^\epsilon_t, \vecc{\mathcal{A}}^\epsilon_t, \mathcal{B}^\epsilon_t)$ in the limit  $\epsilon \to 0$. 

%We will use a perturbative approach to derive, at a formal level, these limiting equations.  The results can be rigorously justified by the weak convergence theorems in \citep{pavliotis2008multiscale}. 

%Moreover, the mode of convergence could be strengthen to a pathwise one and error estimates could be otained along the lines of \citep{Lim2018,birrell2018entropy}.  

%\begin{ass} \label{ass1} The global solutions, defined on %$[0,T]$, to the pre-limit SDEs \eqref{eqx}-\eqref{eqA} and to the limiting SDE \eqref{mainlimitingeqn} a.s. exist and are unique for all $\epsilon > 0$ (i.e. there are no explosions).
%\end{ass}

We now introduce our notation and provide some reminders on transformation of stochastic integrals.  \\

\noindent {\bf Notation.} Consider the diffusion process $\vecc{Z}_t \in \RR^N$, $t \geq 0$,  satisfying the It\^o SDE: 
\begin{equation}
d\vecc{Z}_t = \vecc{b}(t, \vecc{Z}_t) dt + \vecc{\sigma}(t, \vecc{Z}_t) d\vecc{W}_t, 
\end{equation}
where $\vecc{b}: \RR^+ \times \RR^N \to \RR^N$, $\vecc{\sigma}: \RR^+ \times \RR^N  \to \RR^{N \times M}$ (differentiable in $\vecc{Z}$), and $\vecc{W}_t \in \RR^M$ is a Wiener process. Equivalently, it can be cast as the following {\it Stratonovich SDE}: 
\begin{equation}
d\vecc{Z}_t = \vecc{u}(t,\vecc{Z}_t) dt +  \vecc{\sigma}(t, \vecc{Z}_t) \circ d\vecc{W}_t, \label{strat_sde}
\end{equation}
where $\vecc{u}(t,\vecc{Z}_t) = \vecc{b}(t, \vecc{Z}_t)  - \vecc{c}(t,\vecc{Z}_t) $, the symbol  $\circ$ denotes Stratonovich convention (without the symbol $\circ$,  It\^o convention is taken), and, in index-free notation, 
\begin{equation}
\vecc{c} = \frac{1}{2}[ \vecc{\nabla} \cdot (\vecc{\sigma} \vecc{\sigma}^{T}) - \vecc{\sigma} \vecc{\nabla} \cdot ( \vecc{\sigma}^{T}) ].
\end{equation}
In the above, $\vecc{\nabla} \cdot$ denotes divergence operator which contracts a matrix-valued function to a vector-valued function: for the matrix-valued function $\vecc{A}(\vecc{Z})$, the $i$th component of its divergence is given by
\begin{equation}
(\vecc{\nabla} \cdot \vecc{A})^i = \sum_{j=1}^N \frac{\partial A^{ij}}{\partial Z^j}.
\end{equation}
Equivalently, in components,
\begin{equation}
c^i = \frac{1}{2} \frac{\partial \sigma^{ij}}{\partial X^k} \sigma^{kj},
\end{equation}
where $\sigma^{ij}$ denotes the $(i,j)$-entry of the matrix $\vecc{\sigma}$, $Z^k$ the $k$th component of the vector $\vecc{Z}$, and we have used Einstein's summation convention for repeated indices. \\

%In terms of the kinetic or H\"{a}nggi-Klimontovich drift $\hat{\vecc{u}}$, one has:
%\begin{equation}
%d\vecc{Z}_t = \hat{\vecc{u}}(t,\vecc{Z}_t) dt + 2\vecc{c}(t,\vecc{Z}_t) dt + \vecc{\sigma}(t,\vecc{Z}_t) d\vecc{W}_t. 
%\end{equation}

%We also have the following rules for a $L$-dimensional functional defined in terms of the trajectory $(\vecc{Z}_t)_{t\in [0,T]}$: 
%\begin{align}
%\int_0^T \vecc{H}(t,\vecc{Z}_t) \circ d\vecc{Z}_t &= \int_0^T \vecc{H}(t, \vecc{Z}_t) d\vecc{Z}_t +\frac{1}{2} \int_0^T [\vecc{\nabla} \cdot (\vecc{H} \vecc{D}) - \vecc{H} \vecc{\nabla} \cdot \vecc{D} ](t,\vecc{Z}_t) dt,
%\end{align}
%where $\vecc{H} \in \RR^{L \times N}$ and $\vecc{D} := \vecc{\sigma} \vecc{\sigma}^{T} \in \RR^{N \times N}$ is the diffusion matrix. \\

We make the following assumptions on the SDE systems \eqref{eqx}-\eqref{egB}: 

\begin{ass} \label{ass1} The global solutions, defined on $[0,T]$, to the pre-limit SDEs \eqref{eqx}-\eqref{egB} and to the limiting SDEs \eqref{Xlimit}-\eqref{Anomlimit} a.s. exist and are unique for all $\epsilon > 0$ (i.e. there are no explosions).
\end{ass}

\begin{ass} \label{ass2} The matrix-valued functions  $\{\vecc{U}_2(t,\vecc{X}); t \in [0,T], \vecc{X} \in \RR^{n}\}$ are {\it uniformly positive stable}, i.e. all real parts of the eigenvalues of $\vecc{U}_2(t, \vecc{X})$ are bounded from below, uniformly in $t$ and $\vecc{X}$, by a positive constant. 
\end{ass}

\begin{ass} \label{ass3} For $t \in [0,T]$, $\vecc{X} \in \RR^{n}$, and $i=1,2$, the functions $\vecc{u}_i(t,\vecc{X})$,  $\tilde{\vecc{\sigma}}(t,\vecc{X})$, $\vecc{\sigma}(t,\vecc{X})$, $\vecc{r}(t, \vecc{X})$ are continuous and bounded in $t$ and $\vecc{X}$, and Lipschitz in $\vecc{X}$, whereas the functions $\vecc{U}_i(t,\vecc{X})$, $\vecc{P}(t,\vecc{X})$,  $(\vecc{U}_i)_{\vecc{X}}(t,\vecc{X})$, $\vecc{P}_{\vecc{X}}(t,\vecc{X})$ are continuous in $t$, continuously differentiable in $\vecc{X}$, bounded in $t$ and $\vecc{X}$, and Lipschitz in $\vecc{X}$.
Moreover, the  functions $(\vecc{U}_i)_{\vecc{X} \vecc{X}}(t,\vecc{X})$ ($i=1,2$) and $\vecc{P}_{\vecc{X} \vecc{X}}(t,\vecc{X})$ are bounded for every $t \in [0,T]$ and $\vecc{X} \in \RR^{n}$.
\end{ass}

\begin{ass} \label{ass4} The initial condition $\vecc{X}^\epsilon_0 = \vecc{X}^\epsilon \in \RR^{n}$ is an $\mathcal{F}_0$-measurable random variable that may depend on $\epsilon$, and we assume that $\mathbb{E}[|\vecc{X}^\epsilon|^p] = O(1)$ as $\epsilon \to 0$ for all $p>0$. Also, $\vecc{X}^\epsilon$ converges, in the limit as $\epsilon \to 0$, to a random variable $\vecc{X}$ as follows:
$\mathbb{E}\left[|\vecc{X}^\epsilon - \vecc{X}|^p \right]  = O(\epsilon^{p r_0})$, where $r_0 > 1/2$ is a constant, as $\epsilon \to 0$.  The same conditions are assumed for $\vecc{A}^\epsilon_0$.    The initial condition $\vecc{Y}^\epsilon_0 = \vecc{Y}^\epsilon \in \RR^{m}$ is an $\mathcal{F}_0$-measurable random variable that may depend on $\epsilon$, and we assume that for every $p>0$, $\mathbb{E}[ |\epsilon \vecc{Y}^\epsilon|^p] = O(\epsilon^\alpha)$ as $\epsilon \to 0$, for some $\alpha \geq p/2$.
\end{ass}

The following theorem follows from a straightforward application of Theorem \ref{mainthm}. The last statement in the theorem follows from the proof of Theorem \ref{mainthm} (see \citep{2019Lim} for details). 
\begin{thm}  \label{mainthm2}
Under the Assumption \ref{ass1}-\ref{ass4}, in the limit $\epsilon \to 0$,  the family of processes $(\vecc{X}_t^\epsilon, \vecc{\mathcal{A}}_t^\epsilon)$, $t \in [0,T]$,  converges to $(\vecc{X}_t, \vecc{\mathcal{A}}_t)$ solving the It\^o SDE:
\begin{align}
d\vecc{X}_t &= [\vecc{u}_1(t,\vecc{X}_t) + \vecc{U}_1(t,\vecc{X}_t) \vecc{U}_2^{-1}(t,\vecc{X}_t) \vecc{u}_2(t,\vecc{X}_t)] dt + \vecc{S}_{Ito}(t,\vecc{X}_t) dt + \tilde{\vecc{\sigma}}(t,\vecc{X}_t) d\tilde{\vecc{W}}_t \nonumber \\  &\ \ \ \  + \vecc{U}_1(t,\vecc{X}_t) \vecc{U}_2^{-1}(t,\vecc{X}_t) \vecc{\sigma}(t,\vecc{X}_t) d\vecc{W}_t, \label{Xlimit} \\  
d\vecc{\mathcal{A}}_t &= \vecc{r}(t,\vecc{X}_t) dt + \vecc{P}(t,\vecc{X}_t) d\vecc{X}_t + d\vecc{\mathcal{A}}_t', \label{Alimit}\\
d\vecc{\mathcal{A}}_t' &= [\vecc{\nabla} \cdot (\vecc{P}(t,\vecc{X}_t) \vecc{U}_1(t,\vecc{X}_t) \vecc{\mu}(t,\vecc{X}_t) \vecc{U}^T_1(t,\vecc{X}_t)), \nonumber  \\ 
&\ \ \ \ - \vecc{P}(t,\vecc{X}_t) \vecc{\nabla} \cdot (\vecc{U}_1(t,\vecc{X}_t) \vecc{\mu}(t,\vecc{X}_t) \vecc{U}_1^T(t,\vecc{X}_t)) ] dt, \  \text{  or, in component:} \label{Anomlimit} \\
d(\mathcal{A}'_t)^k &= U_1^{ia} U_1^{jb} (U_2^{-1} J)^{ab} \frac{\partial P^{ki}}{\partial X^j} dt.  \label{Anomlimit2} 
\end{align}
In the above, $\vecc{S}_{Ito}$ is the noise-induced drift:
\begin{equation}
\vecc{S}_{Ito} = \vecc{\nabla} \cdot (\vecc{U}_1 \vecc{U}_2^{-1} \vecc{J} \vecc{U}_1^T) - \vecc{U}_1 \vecc{U}_2^{-1} \vecc{\nabla} \cdot (\vecc{J} \vecc{U}_1^T), 
\end{equation}
with $\vecc{J}$ solving the Lyapunov equation
\begin{equation} \label{lyap}
\vecc{U}_2 \vecc{J} + \vecc{J} \vecc{U}_2^T = \vecc{\sigma} \vecc{\sigma}^T, 
\end{equation} 
and $\vecc{\mu} = \vecc{U_2^{-1}} \vecc{J}$. 
The convergence is in the following sense: for all finite $T>0$,
\begin{equation}
 \sup_{t \in [0,T]} |\vecc{X}^\epsilon_t - \vecc{X}_t| \to 0, \ \ \sup_{t \in [0,T]} |\vecc{\mathcal{A}}^\epsilon_t - \vecc{\mathcal{A}}_t| \to 0,  \end{equation}
in probability, in the limit as $\epsilon \to 0$. The family of functionals $\mathcal{B}_t^\epsilon = \frac{\epsilon}{2} |\vecc{Y}_t^\epsilon|^2 $ converges to $Tr(\vecc{J}(t,\vecc{X}_t))$ as $\epsilon \to 0$ in the following sense:  for all finite $T > 0$, 
\begin{equation}
\sup_{t \in [0,T]} \int_0^t | \mathcal{B}_s^\epsilon - Tr(\vecc{J}(s, \vecc{X}_s)| ds \to 0
\end{equation}
in probability as $\epsilon \to 0$. 
\end{thm} 

The following remarks describe the link between {\it symmetry breaking} (violation of a detailed balance condition) and {\it area anomaly} (concerning the appearance of the anomalous contributions, $\vecc{S}_{Ito} dt$ and $d\vecc{\mathcal{A}}'_t$, in the homogenized equations). 

\begin{rmk} \label{rmk41}
We recall some connections to relevant concepts from nonequilibrium statistical mechanics \citep{pavliotis2014stochastic}. Define the matrix $\vecc{\mu}$ and $\vecc{\nu}$, by 
\begin{align}
\mu^{ab} &:= \int_0^\infty \mathbb{E} Y_\tau^a Y_0^b  d\tau, \\ 
2 \mu_S^{ab} &:= \mu^{ab}+ \mu^{ba} =: \nu^{ac} \nu^{bc}.
\end{align}  
Let $\vecc{L}_0$ be the infinitesimal generator corresponding to the fast dynamics in $\vecc{Y}$, i.e. $\vecc{L}_0 = -\vecc{U}_2(t,\vecc{X}) \vecc{Y} \cdot \vecc{\nabla}_{\vecc{Y}} + \frac{1}{2} (\vecc{\sigma}(t,\vecc{X}) \vecc{\sigma}^T(t,\vecc{X})): \vecc{\nabla}_{\vecc{Y}} \vecc{\nabla}_{\vecc{Y}}$, where 
$\vecc{A} : \vecc{\nabla}_{\vecc{Y}} \vecc{\nabla}_{\vecc{Y}} := \sum_{i,j} A^{ij} \frac{\partial^2}{\partial Y^i \partial Y^j}$.   
Using the time integral representation formula for $(-\vecc{L}_0)^{-1}$, one finds  $\mu^{ab}  = \overline{Y^b (-\vecc{L}_0^{-1} Y^a)}$, where overbar denotes averaging with respect to the invariant density of a mean zero Gaussian process with the covariance matrix $\vecc{J}$. This is an example of the Green-Kubo formula, which is important for the calculation of transport coefficients \citep{pavliotis2010asymptotic}.  It is straightforward to compute that $\vecc{\mu} = \vecc{U}_2^{-1} \vecc{J}$ and $\vecc{\nu} = \vecc{U}_2^{-1} \vecc{\sigma}$. Recall that $\vecc{J}$ solves the Lyapunov equation \eqref{lyap}, which can be rewritten as $\vecc{L} + \vecc{L}^T = \vecc{D}$, where $\vecc{L} := \vecc{U}_2 \vecc{J}$ is the Onsager matrix of kinetic coefficient (associated to the fast dynamics) and $\vecc{D} = \vecc{\sigma} \vecc{\sigma}^T$ is the diffusion matrix \citep{godreche2018characterising}. 

It is well known that the {\it detailed balance condition} (the condition for the fast process to be reversible, or equivalently, for its infinitesimal generator to be symmetric), for a given $t$ and $\vecc{X}$, holds if and only if  $\vecc{U}_2 \vecc{D}$ is symmetric, i.e. $\vecc{U}_2 \vecc{D} = \vecc{D} \vecc{U}_2^T$ \citep{gardiner2009stochastic}.  In this case, the stationary covariance matrix is $\vecc{U}_2^{-1} \vecc{D}/2$ and the corresponding stationary state is an equilibrium one. In particular, this symmetry condition implies that $\vecc{\mu}$ is symmetric and $\vecc{\mu} = \vecc{\mu}_S$. The converse is not true unless $\vecc{U}_2^2 \vecc{J}$ is symmetric. When the symmetry condition is broken, the fast process is irreversible and has a nonequilibrium stationary state. One can quantify the irreversibility of the process as follows. We write $\vecc{L} = \vecc{D}/2+\vecc{Q}$ and $\vecc{L}^T = \vecc{D}/2 - \vecc{Q}$ so that we can use $\vecc{Q} = (\vecc{L}-\vecc{L}^T)/2$, {\it the antisymmetric part of the Onsager matrix}, to measure the irreversibility of the fast process. If the fast process is reversible, then the Onsager matrix $\vecc{L} = \vecc{D}/2$ is symmetric and $\vecc{Q}=\vecc{0}$. We refer to \citep{godreche2018characterising,Macieszczak18} and the references therein for a list of works on quantification of the asymmetry of the Onsager matrix. 
\end{rmk}

\begin{rmk} In the case when $\tilde{\vecc{\sigma}} = \tilde{\vecc{\sigma}}(t)$ and $\vecc{\sigma} = \vecc{\sigma}(t)$ are independent of the state, we have:
\begin{align}
d\vecc{X}_t &= (\vecc{u}_1(t,\vecc{X}_t) + \vecc{U}_1(t,\vecc{X}_t) \vecc{U}_2^{-1}(t,\vecc{X}_t) \vecc{u}_2(t,\vecc{X}_t) + d\vecc{X}_t''  \nonumber \\  
&\ \ \  \  + \tilde{\vecc{\sigma}}(t) d \tilde{\vecc{W}}_t + \vecc{U}_1(t,\vecc{X}_t) \vecc{U}_2^{-1}(t,\vecc{X}_t) \vecc{\sigma}(t) \circ d\vecc{W}_t,
\end{align}
with $d\vecc{X}_t'' =  \vecc{H}_{Str}(t,\vecc{X}_t)) dt$, where  $\vecc{H}_{Str}$ is the additional drift term which can be written in two equivalent ways. The first one is in terms of $\vecc{Q}$, $\vecc{L}$ and $\vecc{\nu}$ introduced earlier and  $\vecc{H}_{Str}$ is written compactly as a sum of three contributions:
\begin{equation}
\vecc{H}_{Str} = \vecc{\nabla} \cdot (\vecc{U}_1 \vecc{U}_2^{-1} (\vecc{U}_1 \vecc{U}_2^{-1} \vecc{Q})^T) - \vecc{U}_1 \vecc{U}_2^{-1} \vecc{\nabla}\cdot ((\vecc{U}_1 \vecc{U}_2^{-1} \vecc{L})^T) + \frac{1}{2}(\vecc{U}_1 \vecc{U}_2^{-1}) \vecc{\sigma} \vecc{\nabla} \cdot ((\vecc{U}_1 \vecc{\nu})^T).
\end{equation}
The second way is in terms of $\vecc{Q}$, Lie brackets of vector fields and $\vecc{\nu}$: 
\begin{align}
H_{Str}^i &= \frac{\partial (U_1 U_2^{-1})^{ip}}{\partial X^k} (U_1 U_2^{-1})^{kl} Q^{lp} = \frac{1}{2} Q^{lp}[\vecc{G}_l, \vecc{G}_p]^i,
\end{align}
where the vector fields  $\vecc{G}_l$ are associated to the $l$th column of the matrix $\vecc{U}_1\vecc{U}_2^{-1}$ and $[\cdot, \cdot]$ denotes the Lie bracket of two vector fields. The {\it antisymmetric} matrix $\vecc{Q}$ (which, as discussed earlier, measures the irreversibility of the fast process) encodes the stochastic area of the limiting dynamical process, and  $\vecc{H}_{Str}$ would vanish in the one-dimensional case (c.f. \citep{ikeda2014stochastic}, or Section 2 in \citep{lejay2003importance} for the point of view of interpolation problem for trajectories). 

The irreversibility of the fast process generates macroscopic current in the stationary state and induces some loops in the trajectories. It turns out that the area generated by these loops is of $O(1)$ as $\epsilon \to 0$. As a result, zooming in the small scale $\vecc{X}_t$ ``spins'' around a modified mean trajectory \citep{lejay2002convergence, lejay2003importance}. We refer the reader to Section \ref{app_levy} for an illustration of such phenomenon in a simple example. The phenomena of area anomaly has  been discovered and studied recently in different problem settings  \citep{chevyrev2016multiscale,lopusanschi2017area, lopusanschi2018levy} (see also the references therein). One rigorous framework for understanding these phenomena is based on the theory of rough paths \citep{lyons1998differential,friz2010multidimensional,friz2020course}.

%We remark that in the case when $\vecc{\sigma}$ is state-dependent, we see that even when the fast process admits an equilibrium stationary state  or in the one-dimensional setting (i.e. $\vecc{Q} = \vecc{0}$), $\vecc{H}$ does not necessarily vanish, in which case the limiting equation is not in the Stratonovich form. Therefore,  state dependence of $\vecc{\sigma}$ would give another contribution to the modified drift term.
\end{rmk}

\begin{rmk}
The evolution of the effective functional is described by:
\begin{align}
d\vecc{\mathcal{A}}_t &= \vecc{r} dt + \vecc{P}  \circ d\vecc{X}_t + d\vecc{\mathcal{A}}_t'', \\ 
d\vecc{\mathcal{A}}''_t &=  \left[\vecc{\nabla} \cdot \left(\vecc{P} \left(\vecc{U}_1 \vecc{\mu}_A^T  \vecc{U}_1^T - \frac{1}{2} \tilde{\vecc{\sigma}} \tilde{\vecc{\sigma}}^T \right) \right) - \vecc{P} \vecc{\nabla} \cdot \left(\vecc{U}_1 \vecc{\mu}_A^T \vecc{U}_1^T  - \frac{1}{2} \tilde{\vecc{\sigma}} \tilde{\vecc{\sigma}}^T \right) \right]  dt,
\end{align}
where $\vecc{\mu}_A$ is the antisymmetric part of $\vecc{\mu}$. In component form, we have:
\begin{align}
d(\mathcal{A}_t'')^i &= \frac{1}{2} U_1^{kb} U_1^{ja} \mu_A^{ab} \left( \frac{\partial P^{ij}}{\partial X^k} - \frac{\partial P^{ik}}{\partial X^j} \right) dt - \frac{\partial P^{ij}}{\partial X^k} (\tilde{\sigma} \tilde{\sigma}^T)^{kj} dt.  
\end{align}

Therefore,  whenever $\vecc{\mu}_A = \vecc{0}$ (a sufficient condition for this is when $\vecc{Q} = \vecc{0}$ and $\tilde{\vecc{\sigma}} = \vecc{0}$), $d\vecc{\mathcal{A}}_t'' = \vecc{0}$ and the effective SDE for the functional $\vecc{\mathcal{A}}_t$ can be expressed entirely in terms of the trajectory of the slow process in the {\it Stratonovich prescription}. Otherwise, the loops induced by irreversibility of the fast dynamics in the $\vecc{X}$-trajectory generally cause $\vecc{\mathcal{A}}_t$,  a functional of the $\vecc{X}$-trajectory, to ``spin'' around a modified mean trajectory in the limit. Similar results, albeit in a different and more abstract context, were also shown and discussed in \citep{lejay2003importance}. In the very special case when $\vecc{r} = \vecc{0}$, $\vecc{P}$ is an identity matrix and $\vecc{\mathcal{A}}_0^\epsilon = \vecc{X}_0^\epsilon = \vecc{0}$, we have $\vecc{\mathcal{A}}_t^\epsilon = \vecc{X}_t^\epsilon$ and therefore the effective description for both dynamical variable and functional coincides -- see Remark \ref{rmk41} for expression of the anomalous contribution in this case. 

Finally, we remark that even in the general case when $\vecc{\mu}_A$ is non-zero, the effective SDE for the functional $\mathcal{\vecc{A}}_t$ can be expressed entirely in terms of the trajectory of the slow process (albeit generally not in the Stratonovich prescription), and therefore the area anomaly due to $\vecc{\mathcal{A}}_t$  here is different from the entropy anomaly studied in \citep{bo2017multiple}, where new independent noise terms need to be introduced in the effective equation for the entropy production.
\end{rmk}

%Lastly, we mention the connection with the Levy area correction. Let consider the case when $n=m=2$, $u_0^i(x)=0$, $U_1(x) = I$, where $I$ is an identity matrix, $\sigma_0 = 0$, $B_2 = \Gamma$, $\sigma_1 = \sigma$, $b_1 = 0$, $r=0$, and $p(x,y)=(-x^2, x^1)$, i.e. when the functional under consideration is $S_t = \int_0^t (x_s^1 dx_s^2 - x_s^2 dx_s^1)$ and the $x_t^i$ are approximants of a scaled Wiener process as $\epsilon \to 0$. 

%One can generalize the above derivation to study the limiting behavior of the antisymmetric matrix-valued functional $A_t^{ij} = \int_0^t (x_s^i dx_s^j - x_s^j dx_s^i)$ (recall that so far we have considered only scalar-valued functionals), where the $x_t$ are approximants of a scaled Wiener process satisfying $dx_t^i = \frac{1}{\epsilon} y_t^i dt$, $dy_t^i = -\frac{1}{\epsilon^2} \Gamma^{ij} y_t^j dt + \frac{1}{\epsilon} \sigma^{ia} dW_t^a$. 

%Mathematically, Levy area correction and symmetry (detailed balance).  Physically, Green-Kubo formula for transport coefficient. 

\section{An Example: Stochastic Area as Work Functional  and Its Homogenization} \label{app_levy}

In this section, we are going to discuss the  area anomaly phenomenon  and its consequences in the context of simple Langevin systems describing the motion of a Brownian particle in magnetic field. In the next section, we will study how such phenomena manifests itself for functionals along trajectories of a wide class of multi-dimensional generalized Langevin systems approximating, in various time scale separation scenarios, that of an effective Langevin system. 

Let $(\vecc{q}_s)_{s \geq 0}$, $\vecc{q}_s=(q_s^1, q_s^2) \in \RR^2$, be a stochastic process. Let $\vecc{q}_0 = \vecc{0}$. The {\it stochastic area of  $(\vecc{q}_s)_{s \in [0,t]}$  on the interval $[0,t]$} is defined as the random variable:
\begin{equation}
S_t = \frac{1}{2} \int_0^t (q_s^1 dq_s^2 - q_s^2 dq_s^1). 
\end{equation} 
Viewing $t$ as a continuous-time parameter, this gives rise to the area process $(S_t)_{t \geq 0}$. The above formula, with $\vecc{q}_s = \vecc{W}_s$ (i.e. a 2D Wiener process), is an object first introduced and studied by L\'evy in \citep{levy1951wiener}. His formula formally defines the area (which is random) included by the curve $\mathcal{C}_t=\{Q^1=q^1_s, Q^2 = q^2_s, \ s \in [0, t]\}$ and its chord. Extension of this definition to the case when $\vecc{q}_s \in \RR^d$, $d>2$, is straightforward \citep{malham2010introduction}.

Let $(\vecc{\eta}^\epsilon_s)_{s \in [0, T]}$, $\epsilon >0$, be a family of sufficiently smooth approximations of the Wiener process $(\vecc{W}_s)_{s \in [0, T]}$, where $\vecc{\eta}_s^\epsilon$ converges to $\vecc{W}_s$ as $\epsilon \to 0$  in a pathwise sense. A natural question is then whether or not the stochastic area of $(\vecc{\eta}_s)_{s \in [0,t]}$ converges to L\'evy's stochastic area as $\epsilon \to 0$. We will show that this is generally not true and discuss the consequences in the context of a physical system. We would expect similar conclusion to hold had we replaced $S_t$ with other functionals. 

Consider the motion of a charged (non-relativistic) particle undergoing Brownian motion in the presence of a magnetic field. Such motion is of interest in astrophysics, as motion from interacting charged particles produces observed light curves with interesting peculiarities \citep{Harko2016}. For simplicity, here we consider the case where the magnetic field, $\vecc{B}$, points in the $z$-direction with a constant magnitude $B$ and study the motion of the particle in the 2D plane perpendicular to the magnetic field\footnote{The analysis beyond this case  is straightforward but involves richer physics. For instance, the charged particle may spiral in a non-trivial configuration-dependent manner when the magnetic field is position-dependent and points in arbitrary direction.}. In the absence of external forces and noise where the magnetic field is the dominant factor determining the motion, the particle revolves in a circular orbit with a  frequency $\Omega$, producing  current loops. In this case, the magnetic force is $\vecc{F}_B = \Omega \vecc{V} \times \vecc{e}_3 = \Omega(v^2, -v^1,0)$, where $\vecc{V}=(v^1, v^2,v^3) \in \RR^3$ is the velocity of the charged particle, $\times$ denotes  cross product  and $\vecc{e}_3=(0,0,1)$. It does no work on the particle, even though the direction of motion of the particle is changed. 

Taking into account this magnetic force, as well as a drag and noise term to model collisions of the charged particle with surrounding particles, the evolution of  position $\vecc{q}_t = (q^1_t,q^2_t)$ and velocity $\vecc{v}_t = (v^1_t, v^2_t)$  of the particle on the 2D plane (assuming $\vecc{q}_0 = \vecc{0}$) can be described by the SDE:
\begin{align}
d\vecc{q}_t &= \vecc{v}_t dt, \label{qq} \\ 
m d\vecc{v}_t &= -\Omega \vecc{J}\vecc{v}_t dt - \vecc{v}_t dt + \vecc{A} d\vecc{W}_t,   \label{pp}
\end{align}
where $m>0$ is the mass of the particle,  $\Omega =\frac{qB}{c}$ (with $q$ the charge of the particle, $c$ the speed of light and $B$ the magnitude of the constant magnetic force) is the Lamor frequency (up to a multiplicative factor of $1/m$), $\vecc{A} = \vecc{I}+\Omega \vecc{J}$ $\bigg($with $\vecc{I}$ identity matrix and $\vecc{J} = \begin{bmatrix} 0 & -1 \\ 1 & 0 \end{bmatrix} \bigg)$, and $\vecc{W}_t$ is a Wiener process.  Note that $\vecc{A}$ is positive stable (but not symmetric unless $\Omega = 0$) and $ \vecc{A}\vecc{W}_t$ is a Brownian motion with the covariance matrix $(1+\Omega^2)\vecc{I}$.

Let us now suppose that the charged particle is additionally subject to an external, non-conservative force field, $\vecc{f}_{nc}(t,\vecc{q})$, so that the equations of motion become:
\begin{align}
d\vecc{q}_t &= \vecc{v}_t dt, \label{qqq} \\ 
m d\vecc{v}_t &= -\Omega \vecc{J} \vecc{v}_t dt - \vecc{v}_t dt + \vecc{f}_{nc}(t,\vecc{q}_t) dt +   \vecc{A} d\vecc{W}_t.  \label{ppp}
\end{align}
In this case,  following the approach in stochastic energetics \citep{sekimoto2010stochastic}, we write the kinetic energy of the charged particle as $\mathcal{E}_t := \frac{1}{2} m v_t^2 = \mathcal{Q}_t + \mathcal{W}_t$, where the heat $\mathcal{Q}_t$ and work $\mathcal{W}_t$ satisfies:
\begin{align}
d\mathcal{Q}_t &= m \vecc{v}_t \circ d\vecc{v}_t - \vecc{f}_{nc}(t,\vecc{q}_t) \cdot d\vecc{q}_t, \\ 
d\mathcal{W}_t &= \vecc{f}_{nc}(t,\vecc{q}_t) \cdot d\vecc{q}_t,
\end{align} 
where $\circ$ denotes Stratonovich integration and $\cdot$ denotes inner product.  {\it In the special case where $\vecc{f}_{nc}(t,\vecc{q}) = \frac{1}{2}(-q^2, q^1)$, the resulting work is exactly the  stochastic area of the position process, i.e. $d\mathcal{W}_t = \frac{1}{2} (\vecc{J} \vecc{q}_t)^T d\vecc{q}_t = dS_t$.} We will work with this special case in the following.

Setting $m = \epsilon$, we now consider the following rescaled family of the system \eqref{qq}-\eqref{pp}, together with the  SDEs defining the stochastic areas, $S_t^\epsilon$, of the (rescaled) position process of the charged particle: 
\begin{align}
d\vecc{q}^\epsilon_t &= \vecc{v}_t^\epsilon dt, \\
\epsilon d\vecc{v}_t^\epsilon &= -\vecc{A} \vecc{v}_t^\epsilon dt+  \vecc{A}d\vecc{W}_t, \\
dS_t^\epsilon &= \frac{1}{2}(\vecc{J} \vecc{q}^\epsilon_t)^T d\vecc{q}^\epsilon_t. 
\end{align}
A straightforward application of Theorem \ref{mainthm2} allows us to find out whether the family of stochastic areas of $(\vecc{q}^\epsilon_s)_{s \in [0,t]}$ converges to L\'evy's stochastic area as $\epsilon \to 0$.

\begin{cor}
In the limit $\epsilon \to 0$, the family of processes $(\vecc{q}_t^\epsilon,  S_t^\epsilon)$ converges to $(\vecc{W}_t, \bar{S}_t)$, where
\begin{equation}
\bar{S}_t = S_t^{Levy}-  \frac{\Omega}{2} t, 
\end{equation}
with $S_t^{Levy}$  L\'evy's stochastic area.
More precisely, for all finite $T > 0$,  $\sup_{t \in [0,T]} |\vecc{q}_t^\epsilon - \vecc{W}_t| $, $\sup_{t \in [0,T]} |S_t^\epsilon - \bar{S}_t| \to 0$ in probability\footnote{In fact, a stronger $L^p$ convergence result (for $p > 0$) can be obtained in this case.}, as $\epsilon \to 0$. 
\end{cor}

Therefore, unless $\Omega = 0$ the stochastic area (which here carries the meaning of work) of the pre-limit process  does not converge to L\'evy's area in the small mass limit, even though the pre-limit process converges to a Wiener process. The correct limiting area (work) includes an additional term (which we refer to as area anomaly) that depends on the frequency at which the charged particle circles around the magnetic field, retaining  in the limit the information on how the charged particle is moving under  presence of the magnetic force. 

The frequency $\Omega$ can be interpreted as a symmetry breaking parameter. Indeed, when $\Omega > 0$, $\vecc{A}$ is not symmetric, and so the irreversibility (breaking of detailed balance) of the fast velocity process generates macroscopic current  in the stationary state and induces loops in the position space whose areas are of $O(1)$ as $\epsilon \to 0$. This irreversibility can be quantified using the antisymmetric part of the Onsager matrix \citep{godreche2018characterising}, which in this case can be computed to be $\vecc{Q} = \frac{1+\Omega^2}{2} \left(\frac{\vecc{A}-\vecc{A}^T}{2}\right) = \Omega \vecc{J}$, whose off-diagonal entries encode the area anomaly.  

From a physical point of view, such phenomenon  may be experimentally realized, along the line of \citep{argun2017experimental}, in a microscopic heat engine generating a torque via circular motion, from which work may possibly be extracted.  On the other hand, rich mathematical insights on  the phenomenon can be obtained using the theory of rough paths \citep{friz2015physical,bruned2016examples}.

\section{Homogenization of GLEs and Their Functionals}  \label{sect_limit}

%Put $\tilde{x}^\epsilon$ etc for variables used for weak convergence results. 

In this section we explore five homogenization procedures for the GLEs and the associated functionals of interest: 
\begin{itemize}
\item[(5.1)] a Markovian limit;
\item[(5.2)] a limit where the small mass limit is taken after the Markovian limit in (5.1);
\item[(5.3)] the small mass limit;  
\item[(5.4)] a limit where a Markovian limit is taken after the small mass limit in (5.3); and
\item[(5.5)] a joint Markovian and small mass limit.
\end{itemize}
For each procedure, we first state the problem,  motivation as well as the assumptions, and then present the results. These results are obtained by applying Theorem \ref{mainthm2}, upon verifying the assumptions. Since the verification is straightforward  we omit the proof for these results.  We then discuss the commutativity of these procedures and the consequences of imposing/breaking various symmetry conditions (including fluctuation-dissipation relations).

For all these homogenization procedures, we are studying the case where the colored noise comes from two independent sources evolving on different time scales. The noise is modeled by $\vecc{\sigma}(\vecc{x}_t) \vecc{\xi}_t =  \vecc{\sigma}(\vecc{x}_t) \vecc{C}_2 \vecc{\beta}_t =  \vecc{\sigma}_s(\vecc{x}_t) \vecc{\xi}_t^{(s)} + \vecc{\sigma}_f(\vecc{x}_t) \vecc{\xi}_t^{(f)}$, where $\vecc{\xi}_t^{(s)} =  \vecc{C}_s  \vecc{\beta}^{(s)}_t$ and $\vecc{\xi}_t^{(f)} = \vecc{C}_f \vecc{\beta}^{(f)}_t$, with $\vecc{\beta}_t^{(s)}$ and $\vecc{\beta}_t^{(f)}$ satisfying SDEs of the form \eqref{gle4} with different damping and diffusion constants, i.e.:
\begin{align}
d\vecc{\beta}_t^{(s)} &= -\vecc{\Gamma}_s \vecc{\beta}_t^{(s)} dt + \vecc{\Sigma}_s d\vecc{W}_t^{(d_s)}, \\ 
d\vecc{\beta}_t^{(f)} &= -\vecc{\Gamma}_f \vecc{\beta}_t^{(f)} dt + \vecc{\Sigma}_f d\vecc{W}_t^{(d_f)}. 
\end{align}
Here $\vecc{\sigma}_f$ is a non-zero matrix, $\vecc{\sigma}_s$ is a possibly zero matrix, $\vecc{W}_t^{(d_s)} \in \RR^{d_s}$ and $\vecc{W}_t^{(d_f)} \in \RR^{d_f}$ are independent Wiener processes on a filtered probability space $(\Omega, \mathcal{F}, \mathcal{F}_t, \mathbb{P})$ satisfying the usual conditions,  and $\vecc{\xi}_t^{(f)}$ denotes the part of the noise whose correlation times are much smaller than those of $\vecc{\xi}_t^{(s)}$ (the superscript $(s)$ and $(f)$ indicate ``slow'' and ``fast'' respectively). The matrices $\vecc{\Gamma}_i$ $(i=s,f)$ are positive stable and $\vecc{M}_i = \vecc{M}_i^T > 0$ satisfies the Lyapunov equation $\vecc{\Gamma}_i \vecc{M}_i + \vecc{M}_i \vecc{\Gamma}_i^T = \vecc{\Sigma}_i \vecc{\Sigma}_i^T$. We denote the covariance of $\vecc{\xi}_t^{(s)}$ and $\vecc{\xi}_t^{(f)}$ as $\vecc{R}_s(t)$ and $\vecc{R}_f(t)$ respectively.

In the case when $\vecc{\sigma}_s$ is zero, $\vecc{\sigma}(\vecc{x}_t) \vecc{\xi}_t = \vecc{\sigma}_f(\vecc{x}_t) \vecc{\xi}_t^{(f)}$, in which case all the noise correlation time scales are small, so taking these time scales to zero performs the full white noise limit for the GLE. Otherwise, not all noise correlation time scales are small and therefore not all of these time scales will be taken to zero, performing only a partial white noise limit for the GLE -- this retains the influence of the colored noise on the system in the limit. Retaining some effects of non-Markovianity when taking limits is important to obtain  more realistic effective models, particularly for anomalously diffusing systems (see \citep{2019Lim}).

We assume, throughout the rest of the paper, that:

\begin{ass} \label{ass_k} The matrices 
\begin{equation}
\vecc{K}_i = \vecc{C}_i \vecc{\Gamma}_i^{-1}\vecc{M}_i \vecc{C}_i^T \ \ (i=1,f)
\end{equation}
 are non-zero and invertible (but not necessarily positive definite). 
\end{ass}
 
This assumption is necessary for a meaningful Markovian limit and implies that the GLE models normal diffusion (see \citep{2019Lim} for  cases where the assumption is violated). The matrix $\vecc{K}_1$ is the effective damping constant and $\vecc{K}_f$ the effective diffusion constant (for the fast noise process $\vecc{\xi}_t^{(f)}$) in the GLE \citep{Lim2018}.

In all cases, we are assuming that there are no explosions, i.e. almost surely, for every $\epsilon > 0$ there exists global unique solution to the pre-limit SDE  system  and also to the limiting SDE system  on the time interval $[0,T]$. Other assumptions needed concern the initial conditions as well as the regularity and boundedness of the coefficients in the GLE. Note that we have chosen to work with a rather strong assumptions here -- they can be relaxed in various directions at an increased cost of technicality but we choose not to pursue this here.  %go and look at the AHP paper

%decide if we want to include time dependence in the coeffs of GLE...

\begin{ass} \label{ass_regb} {\it Regularity and boundedness.}
For $t \in \RR^+$,  $\vecc{y} \in \RR^{d}$, the functions  $\vecc{F}(t, \vecc{y})$, $\vecc{\sigma}_0(\vecc{y})$, $\vecc{\sigma}_s(\vecc{y})$ and $\vecc{\sigma}_f(\vecc{y})$ are continuous and bounded (in $t$ and $\vecc{y}$) as well as Lipschitz in $\vecc{y}$, whereas the functions $\vecc{\gamma}_0(\vecc{y})$, $\vecc{g}(\vecc{y})$, $\vecc{h}(\vecc{y})$, $(\vecc{\gamma}_0)_{\vecc{y}}( \vecc{y})$, $(\vecc{g})_{\vecc{y}}(\vecc{y})$ and $(\vecc{h})_{\vecc{y}}(\vecc{y})$ are continuously differentiable and Lipschitz in $\vecc{y}$ as well as bounded (in $\vecc{y}$).
Moreover, the  functions $(\vecc{\gamma}_0)_{\vecc{y}\vecc{y}}(\vecc{y})$, $(\vecc{g})_{\vecc{y}\vecc{y}}( \vecc{y})$ and $(\vecc{h})_{\vecc{y}\vecc{y}}(\vecc{y})$  are bounded for every $\vecc{y} \in \RR^{d}$. 
\end{ass}

\begin{ass} \label{ass_initcond}
{\it Initial conditions.} The initial data $\vecc{x}, \vecc{v} \in \RR^d$ are $\mathcal{F}_0$-measurable random variables independent of the $\sigma$-algebra generated by the Wiener processes $\vecc{W}^{(d_s)}$ and $\vecc{W}^{(d_f)}$. They are independent of $\epsilon$ and have finite moments of all orders.
\end{ass} 

We may also need one of the following stability assumptions when studying certain procedures:

\begin{ass} \label{ass_gamma0} The matrix-valued functions $\{\vecc{\gamma}_0(\vecc{x}); \vecc{x} \in \RR^d \}$ are uniformly positive stable. 
\end{ass}

\begin{ass} \label{ass_gamma1}
The matrix-valued functions $\{\vecc{\Gamma}(\vecc{x}) = \vecc{\gamma}_0(\vecc{x}) +  \vecc{g}(\vecc{x}) \vecc{K}_1 \vecc{h}(\vecc{x}); \vecc{x} \in \RR^d\}$ are uniformly positive stable. 
\end{ass}

\subsection{A Markovian Limit} \label{sect_limit1}
We introduce the scaling $\vecc{\kappa}(t) \mapsto \frac{1}{\epsilon} \vecc{\kappa}\left(\frac{t}{\epsilon} \right)$ and $\vecc{R}_f(t) \mapsto \frac{1}{\epsilon}  \vecc{R}_f\left(\frac{t}{\epsilon} \right)$, where $\epsilon > 0$ is a small parameter, in the GLE \eqref{genle_general} with the colored noise term $\vecc{\sigma}(\vecc{x}_t) \vecc{\xi}_t = \vecc{\sigma}_s(\vecc{x}_t) \vecc{\xi}_t^{(s)} + \vecc{\sigma}_f(\vecc{x}_t) \vecc{\xi}_t^{(f)}$ as defined above. This is the limit where all the memory time scales, associated with the history-dependent damping term, and the relevant noise correlation time scales tend to zero at the same rate, and is therefore a partial Markovian limit. Our goal is to study the limit $\epsilon \to 0$ of the resulting generalized Langevin dynamics as well as of the work-like and heat-like functional of the system. 

Implementing this scaling, and introducing the auxiliary process 
\begin{equation} 
\vecc{y}^\epsilon_t =  \int_0^t e^{-\vecc{\Gamma}_1(t-s)} \vecc{M}_1 \vecc{C}_1^T \vecc{h}(\vecc{x}^\epsilon_s) \vecc{v}^\epsilon_s ds,
\end{equation}
the process $(\vecc{x}^\epsilon_t, \vecc{v}^\epsilon_t, \vecc{y}^\epsilon_t, \vecc{\beta}^{(f)\epsilon}_t, \vecc{\beta}_t^{(s)\epsilon})$ satisfies the  SDE system: 
\begin{align}
d\vecc{x}^\epsilon_t  &= \vecc{v}^\epsilon_t dt, \label{mark_gle1} \\ 
m d\vecc{v}^\epsilon_t &= \vecc{F}(t,\vecc{x}^\epsilon_t) dt - \vecc{\gamma}_0(\vecc{x}^\epsilon_t) \vecc{v}^\epsilon_t dt -  \vecc{g}(\vecc{x}^\epsilon_t) \vecc{C}_1 \vecc{y}^\epsilon_t dt +  \vecc{\sigma}_{f}(\vecc{x}^\epsilon_t) \vecc{C}_f \vecc{\beta}^{(f)\epsilon}_t dt + \vecc{\sigma}_s(\vecc{x}^\epsilon_t) \vecc{C}_s \vecc{\beta}_t^{(s)\epsilon} dt, \\
\epsilon d\vecc{y}^\epsilon_t &= -\vecc{\Gamma}_1 \vecc{y}^\epsilon_t dt + \vecc{M}_1 \vecc{C}_1^T \vecc{h}(\vecc{x}^\epsilon_t) \vecc{v}^\epsilon_t dt, \\ 
\epsilon d\vecc{\beta}^{(f)\epsilon}_t &= -\vecc{\Gamma}_f  \vecc{\beta}_t^{(f)\epsilon} dt + \vecc{\Sigma}_f d\vecc{W}_t^{(d_f)},  \label{mark_gle4} \\
d\vecc{\beta}^{(s)\epsilon}_t &= -\vecc{\Gamma}_s \vecc{\beta}_t^{(s)\epsilon} dt + \vecc{\Sigma}_s d\vecc{W}_t^{(d_s)}.  \label{mark_gle5}
\end{align}

The heat-like functional $\mathcal{Q}_t$ and work-like functional $\mathcal{W}_t$ satisfy the following  SDEs:
\begin{align}
d\mathcal{Q}^\epsilon_t &= m \vecc{v}^\epsilon_t \cdot d\vecc{v}^\epsilon_t - \vecc{F}(t,\vecc{x}^\epsilon_t) \cdot d\vecc{x}^\epsilon_t , \\ 
d\mathcal{W}^\epsilon_t &= \frac{\partial U}{\partial t} dt + \vecc{f}_{nc}(t,\vecc{x}^\epsilon_t) \cdot d\vecc{x}^\epsilon_t, 
\end{align}
where $(\vecc{x}^\epsilon_t,\vecc{v}^\epsilon_t)$ solves the SDE system \eqref{mark_gle1}-\eqref{mark_gle5}. Note that in the special case of $d=2$ with $U:=0$, $\vecc{f}_{nc}(t,\vecc{x}) := \frac{1}{2} (-x^2, x^1)$, the  work-like functional is simply  stochastic area of the position process and the heat-like functional is the difference between the kinetic energy and this area.

The dynamics in $\vecc{y}^\epsilon$ and $\vecc{\beta}^{(f)\epsilon}$ are an order of magnitude faster than those in $\vecc{x}^\epsilon$, $\vecc{v}^\epsilon$, $\vecc{\beta}^{(s)\epsilon}$, $\mathcal{Q}^\epsilon$ and $\mathcal{W}^\epsilon$, and  one has the following results.

\begin{cor} \label{mark_limit_res1}
Under appropriate assumptions on the initial conditions and the coefficients (i.e. Assumption \ref{ass_k}-\ref{ass_initcond}) of the pre-limit SDEs \eqref{mark_gle1}-\eqref{mark_gle5}, the family of processes $(\vecc{x}^\epsilon_t, \vecc{v}^\epsilon_t, \vecc{\beta}_t^{(s)\epsilon})$, satisfying the SDEs \eqref{mark_gle1}-\eqref{mark_gle5}, converges, as $\epsilon \to 0$,  to  the solution $(\vecc{x}_t, \vecc{v}_t, \vecc{\beta}_t^{(s)})$ of the It\^o SDE system: 
\begin{align}
d\vecc{x}_t &= \vecc{v}_t dt, \label{markovlimit_eqx} \\
m d\vecc{v}_t &= \vecc{F}(t,\vecc{x}_t) dt - \vecc{\Gamma}(\vecc{x}_t) \vecc{v}_t dt + \vecc{\Sigma}(\vecc{x}_t) d\vecc{W}_t^{(d_f)} + \vecc{\sigma}_s(\vecc{x}_t) \vecc{C}_s \vecc{\beta}_t^{(s)} dt,  \label{markovlimit_eqv} \\
d\vecc{\beta}_t^{(s)} &= - \vecc{\Gamma}_s \vecc{\beta}_t^{(s)} dt + \vecc{\Sigma}_s d\vecc{W}_t^{(d_s)},  \label{markovlimit_eqbeta}
\end{align}
where $\vecc{\Gamma} = \vecc{\gamma}_0 +  \vecc{g} \vecc{K}_1 \vecc{h}$ and $\vecc{\Sigma} = \vecc{\sigma}_f \vecc{C}_f \vecc{\Gamma}_f^{-1} \vecc{\Sigma}_f$.  The convergence is in the strong pathwise sense as before. 
\end{cor} 

Note that $\vecc{\Sigma}(\vecc{x}_t) \vecc{W}_t^{(d_f)} = \vecc{\sigma}_f(\vecc{x}_t) \vecc{B}_t$, where $\vecc{B}_t$ is a Brownian motion with covariance $\vecc{K}_f+\vecc{K}_f^T$. Corollary \ref{mark_limit_res1} conveys the intuitive idea that in the (partial) Markovian limit the GLE can be effectively replaced by a Langevin equation without any memory term and is driven by a relevant colored noise process. When the FDR of the second kind (i.e., Relation \ref{ass_fdr}(b)) holds, we have $\vecc{\Gamma} = \vecc{\gamma}_0 + \vecc{\sigma} \vecc{K}_1 \vecc{\sigma}^T$ and $\vecc{\Sigma} = \vecc{\sigma}_1 \vecc{C}_1 \vecc{\Gamma}_1^{-1} \vecc{\Sigma}_1$ in Eq. \eqref{markovlimit_eqx}-\eqref{markovlimit_eqbeta}. Note that $\vecc{\Sigma} \vecc{\Sigma}^T \neq \vecc{\Gamma} + \vecc{\Gamma}^T$ unless $\vecc{\gamma}_0 = \vecc{0}$ and $\vecc{g} = \vecc{h}^T = \vecc{\sigma} = \vecc{\sigma}_1 = \vecc{\sigma}_f$.

\begin{cor} \label{mark_limit_res}
 Let $\vecc{\Theta}_A$ denote the antisymmetric part of the matrix $\vecc{\Theta} = \vecc{\sigma}_f \vecc{K}_f^T \vecc{\sigma}_f^T$, with $\vecc{K}_f = \vecc{C}_f \vecc{\Gamma}_f^{-1} \vecc{M}_f \vecc{C}_f^T$,  where $\vecc{M}_f$ solves the Lyapunov equation  $\vecc{\Gamma}_f \vecc{M}_f + \vecc{M}_f \vecc{\Gamma}_f^T = \vecc{\Sigma}_f \vecc{\Sigma}_f^T$.  Under the same assumptions as in Corollary \ref{mark_limit_res1}, the family of processes $(\mathcal{W}^\epsilon_t,\mathcal{Q}^\epsilon_t)$, converges, as $\epsilon \to 0$, to the solution $(\mathcal{Q}_t, \mathcal{W}_t)$ of the SDEs:
\begin{align}
d\mathcal{Q}_t &=  m \vecc{v}_t \circ d\vecc{v}_t -  \vecc{F}(t, \vecc{x}_t) \cdot d\vecc{x}_t   + d\mathcal{Q}_t^{anom},  \\
d\mathcal{W}_t &= \frac{\partial U}{\partial t} dt + \vecc{f}_{nc}(t,\vecc{x}_t) \cdot d\vecc{x}_t,  \label{markovlimit_eqW}
\end{align}
where 
\begin{equation}
d\mathcal{Q}_t^{anom} =  \frac{1}{m} \vecc{\nabla}_{\vecc{v}} \cdot (\vecc{v}_t^{T} \vecc{\Theta}_A(\vecc{x}_t)) dt, 
\end{equation}
and $(\vecc{x}_t, \vecc{v}_t)$ solves the SDE system \eqref{markovlimit_eqx}-\eqref{markovlimit_eqbeta}.  The convergence is in the strong pathwise sense as before. 
\end{cor} 

\begin{cor}
$d\mathcal{Q}_t^{anom} = 0$ if and only if $\vecc{\mu}_f = \vecc{\Gamma}_f^{-1} \vecc{M}_f$ (or equivalently, $\vecc{K}_f$)  is symmetric. In particular, a sufficient condition for $d\mathcal{Q}_t^{anom} = 0$ is when the fast process $\vecc{\beta}^{(f)}_t$ satisfies the detailed balance condition. 
\end{cor}

Therefore, even if the FDR of the second kind holds, $d\mathcal{Q}_t^{anom}$ needs not be zero. In fact, note that $\vecc{\Theta} = \vecc{\sigma}_f \vecc{C}_f \vecc{M}_f \vecc{\Gamma}_f^{-T} (\vecc{\sigma}_f \vecc{C}_f)^T = \vecc{\sigma}_f \vecc{K}_f^T \vecc{\sigma}_f^T$, which can be related to the Onsager matrix associated to the fast dynamics induced by the noise process.   It can be shown that the matrix $\vecc{\Theta}$ is, at least in the case when $\vecc{\sigma}_f$ is a non-zero constant, the time integral of the correlation function of the stationary colored noise process $\tilde{\vecc{\xi}}_t := \vecc{\sigma}_f \vecc{C}_f \vecc{\beta}_t^{(f)}$, i.e. $\Theta^{ab} = \int_0^\infty E[\tilde{\xi}^a_t \tilde{\xi}^b_0] dt$, which is in general not symmetric. From Corollary \ref{mark_limit_res}, we see that, unless $\vecc{\Theta}_A$ vanishes (i.e. when we are in the one-dimensional setting, or in the multi-dimensional setting with all the matrix-valued coefficients diagonal, or when the fast colored noise process admits an equilibrium stationary state), the effective evolution of the functional $\mathcal{Q}_t$ cannot be expressed solely as a Stratonovich integral over the effective trajectory. Interestingly, in the one-dimensional setting, the Stratonovich discretization is justified even if the fluctuation-dissipation relation of the second kind is violated. In the general case, whether $\vecc{\Theta}_A$ vanishes or not is entirely due to the symmetry associated with the fast driving colored process, and, in particular,  is independent of the details of the memory function and the slower driving noise process.

\subsection{The Markovian Limit Followed by the Small Mass Limit} \label{sect_comm}
%Taking small mass, then Markovian
%Taking Markovian, then small mass
%are the result the same and which one agrees with the joint limit
%$\gamma_0$ and $\sigma_0 > 0$ as regularization? 

%REWRITE BELOW taking into account \beta^f, \beta^s and \gamma_0 v... 

We rescale $m \mapsto m_0 \epsilon $, where $m_0 > 0$ is a proportionality constant,  in \eqref{markovlimit_eqx}-\eqref{markovlimit_eqW}. The resulting SDE system then becomes: 
\begin{align}
d\vecc{x}^\epsilon_t &= \vecc{v}^\epsilon_t  dt, \label{res_l3x} \\ 
\epsilon  d\vecc{v}^\epsilon_t &= \vecc{F}(t, \vecc{x}^\epsilon_t) dt - \vecc{\Gamma}(\vecc{x}^\epsilon_t) \vecc{v}^\epsilon_t dt + \vecc{\Sigma}(\vecc{x}^\epsilon_t) d\vecc{W}^{(f)}_t + \vecc{\sigma}_s(\vecc{x}_t^\epsilon) \vecc{C}_s \vecc{\beta}_t^{(s)\epsilon} dt, \label{res_l3u} \\
d\vecc{\beta}_t^{(s)\epsilon} &= - \vecc{\Gamma}_s \vecc{\beta}_t^{(s)\epsilon} dt + \vecc{\Sigma}_s d\vecc{W}_t^{(d_s)}, \label{res_l3beta} \\
d\mathcal{Q}^\epsilon_t &=   m_0 \epsilon \vecc{v}^\epsilon_t \circ d\vecc{v}^\epsilon_t - \vecc{F}(t, \vecc{x}^\epsilon_t) \cdot d\vecc{x}^\epsilon_t + \frac{1}{m_0 \epsilon} \vecc{\nabla}_{\vecc{v}^\epsilon} \cdot ((\vecc{v}^\epsilon_t)^T \vecc{\Theta}_A(\vecc{x}^\epsilon_t)) dt, \label{res_l3Q} \\ 
d\mathcal{W}^\epsilon_t &= \frac{\partial U}{\partial t} dt + \vecc{f}_{nc}(t, \vecc{x}^\epsilon_t) \cdot d\vecc{x}^\epsilon_t.   \label{res_l3W}
\end{align}
We are going to study the limit $\epsilon \to 0$ of the above system. This corresponds to taking the small mass limit after the Markovian limit is taken on the GLE \eqref{genle_general}. We assume that Assumption \ref{ass_gamma1} holds, which is crucial to ensure that the small mass limit of the system described by \eqref{markovlimit_eqx}-\eqref{markovlimit_eqbeta}  is well defined \citep{Lim2018}. 

\begin{cor} \label{follow_limit_res1}
Under appropriate assumptions on the initial conditions and the coefficients (i.e. Assumption \ref{ass_k}-\ref{ass_initcond}) of the pre-limit SDEs \eqref{res_l3x}-\eqref{res_l3u} and Assumption \ref{ass_gamma1}, the family of  processes $\vecc{x}^\epsilon_t$, satisfying the SDEs \eqref{res_l3x}-\eqref{res_l3beta}, converges, as $\epsilon \to 0$, to the solution of the following It\^o SDE: 
\begin{align}
d\vecc{x}_t &= \vecc{\Gamma}^{-1}(\vecc{x}_t) (\vecc{F}(t,\vecc{x}_t) dt  + \vecc{\Sigma}(\vecc{x}_t) d\vecc{W}^{(d_f)}_t + \vecc{\sigma}_s(\vecc{x}_t) \vecc{C}_s \vecc{\beta}_t^{(s)} dt) + \vecc{H}(\vecc{x}_t) dt, \label{difflimit_eqx} \\
d\vecc{\beta}_t^{(s)} &= -\vecc{\Gamma}_s \vecc{\beta}_t^{(s)} dt + \vecc{\Sigma}_s d\vecc{W}_t^{(d_f)}, \label{difflimit_eqbeta} 
\end{align}
where $\vecc{\Gamma} = \vecc{\gamma}_0 + \vecc{g} \vecc{K}_1 \vecc{h}$, $\vecc{\Sigma} = \vecc{\sigma}_f \vecc{C}_f \vecc{\Gamma}_f^{-1} \vecc{\Sigma}_f$,  and $\vecc{H}$ is the noise-induced drift whose expression is given by:
\begin{align}
\vecc{H} &= \vecc{\nabla} \cdot (\vecc{\Gamma}^{-1} \vecc{J}) - \vecc{\Gamma}^{-1} \vecc{\nabla} \cdot \vecc{J}, \label{nid2} 
\end{align}
 where $\vecc{J}$ solves the Lyapunov equation $\vecc{\Gamma} \vecc{J} + \vecc{J} \vecc{\Gamma}^T = \vecc{\Sigma}\vecc{\Sigma}^T = \vecc{\Theta}+\vecc{\Theta}^T$, with 
 
\noindent $\vecc{\Theta} = \vecc{\sigma}_f \vecc{C}_f \vecc{M}_f \vecc{\Gamma}_f^{-T} (\vecc{\sigma}_f\vecc{C}_f)^T$, which was first introduced in Collorary \ref{mark_limit_res}. The convergence is in the strong pathwise sense as before.

If $\vecc{\Gamma} \vecc{\Sigma} \vecc{\Sigma}^T$ is symmetric (detailed balance), then $\vecc{J} = \vecc{\Gamma}^{-1} \vecc{\sigma}_f \vecc{K}_f^{T} \vecc{\sigma}_f^T$ and $\vecc{H}$ simplifies to:
\begin{align}
\vecc{H} &= \vecc{\nabla} \cdot (\vecc{\Gamma}^{-2} \vecc{\sigma}_f \vecc{K}_f^{T} \vecc{\sigma}_f^T ) - \vecc{\Gamma}^{-1} \vecc{\nabla} \cdot (\vecc{\Gamma}^{-1} \vecc{\sigma}_f \vecc{K}_f^{T} \vecc{\sigma}_f^T). \label{nid2_db} 
\end{align}
\end{cor} 

We now discuss the physical meaning of this result. Corollary \ref{follow_limit_res1} says that taking the Markovian limit and then the small mass limit on the GLE \eqref{genle_general} leads to  effective dynamics that can be described by a Langevin equation driven by a white noise and a colored noise $\vecc{\beta}_t^{(s)}$, but with an additional drift correction term $\vecc{H}$ and without a memory term. The additional drift term arises due to the state dependence of $\vecc{\Gamma}$ and needs to be considered in nonequilibrium force measurements   \citep{volpe2016effective}. In other words, the effective dynamics can be equivalently described by:
\begin{align}
\dot{\vecc{x}}_t &= \vecc{\Gamma}^{-1}(\vecc{x}_t) (\vecc{F}(t,\vecc{x}_t)  + \vecc{\Sigma}(\vecc{x}_t) \dot{\vecc{W}}^{(d_f)}_t + \vecc{\sigma}_s(\vecc{x}_t) \vecc{C}_s \vecc{\beta}_t^{(s)} ) + \vecc{H}(\vecc{x}_t), \label{88}
\end{align}
with $\vecc{\Gamma} = \vecc{\gamma}_0 + \vecc{g}\vecc{K}_1 \vecc{g}^T$ when the FDR of the second kind holds. Note that without the colored noise term this is an overdamped Langevin equation with an effective drift term, whose expression can be simplified when $\vecc{\Gamma} \vecc{\Sigma}\vecc{\Sigma}^T$ is symmetric.  Since not all the noise (bath) correlation time scales are taken to zero, the limit retains some non-Markovian effects, giving rise to the appearance of a colored noise in the limit. From modeling point of view, this gives  more realistic effective dynamics compared to those driven by only idealized white noise terms (c.f. \citep{2019Lim}).

%Solving this gives $\vecc{J} = \vecc{\Gamma}^{-1} \vecc{\sigma}_f \vecc{K}_f^{T} \vecc{\sigma}_f^T$ and plugging in this gives the formula for $\vecc{H}$. Note that positive stability of $\vecc{\Gamma}$ ensures that the solution to the Lyapunov equation exists and is unique \citep{bellman1997introduction}. 

\begin{cor} \label{follow_limit_res}
Assume that $\vecc{\Theta}_A = \vecc{0}$ (i.e. $\vecc{\mu}_f = \vecc{\Gamma}_f^{-1} \vecc{M}_f$ is symmetric).   Let $\vecc{K}_A$ denote the antisymmetric part of the matrix $\vecc{K} = \vecc{\Gamma}^{-2} \vecc{\sigma}_f \vecc{K}_f^{T} \vecc{\sigma}_f^T = \vecc{\Gamma}^{-2} \vecc{\Theta}$. Then, under the same assumptions as in Corollary \ref{follow_limit_res1},  as $\epsilon \to 0$, the family of processes $(\mathcal{W}^\epsilon_t, \mathcal{R}^\epsilon_t)$, satisfying the SDEs \eqref{res_l3Q}-\eqref{res_l3W}, converges to the solution of the following SDEs: 
\begin{align}
d\mathcal{W}_t &= \frac{\partial U}{\partial t} dt + \vecc{f}_{nc}(t,\vecc{x}_t) \circ  d\vecc{x}_t + d\mathcal{W}_t', \label{difflimit_eqW}  \\
d\mathcal{R}_t &= \vecc{F}(t, \vecc{x}_t) \circ d\vecc{x}_t + d\mathcal{R}_t',
\end{align}
where 
\begin{align}
d\mathcal{W}_t' &= [\vecc{\nabla} \cdot (\vecc{f}_{nc}^{T}(t,\vecc{x}_t) \vecc{K}_A^{T}(\vecc{x}_t)) -  \vecc{f}_{nc}^{T}(t,\vecc{x}_t)  \vecc{\nabla} \cdot \vecc{K}_A^{T}(\vecc{x}_t)] dt, \\
d\mathcal{R}_t' &= [\vecc{\nabla} \cdot (\vecc{F}^{T}(t,\vecc{x}_t) \vecc{K}_A^{T}(\vecc{x}_t)) -  \vecc{F}^{T}(t,\vecc{x}_t)  \vecc{\nabla} \cdot \vecc{K}_A^{T}(\vecc{x}_t)] dt, 
\end{align}
and $\vecc{x}_t$ solves the SDE \eqref{difflimit_eqx}. The convergence is in the strong pathwise sense as before.
\end{cor} 

Note that if, in addition, $\vecc{\gamma}_0 = \vecc{0}$, then imposing FDR of the second kind guarantees $d\mathcal{W}_t' = d\mathcal{R}_t' = 0$. This is implied by  the following corollary. However,  imposing FDR of the second kind alone does not guarantee $d\mathcal{W}_t' = d\mathcal{R}_t' = 0$ without assuming that $\vecc{\Theta}$ (or $\vecc{\mu}_f$) is symmetric.

\begin{cor} \label{spec1}
Suppose that the assumptions in Corollary \ref{follow_limit_res} holds. Then $d\mathcal{W}_t'=d\mathcal{R}_t = 0$ when $\vecc{\gamma}_0 = \vecc{0}$, $\vecc{g} \propto \vecc{h}^T = \vecc{\sigma}_f$ and $\vecc{K}_f = \vecc{K}_1$. 
\end{cor}

One can write $\vecc{K}$, using the solution $\vecc{J}$ of the Lyapunov equation, explicitly as:
\begin{equation}
\vecc{K} = (\vecc{\gamma}_0 + \vecc{g}\vecc{K}_1 \vecc{h})^{-1} \int_0^\infty e^{-(\vecc{\gamma}_0 + \vecc{g}\vecc{K}_1 \vecc{h}) y} (\vecc{\Theta}+\vecc{\Theta}^T) e^{-(\vecc{\gamma}_0 + \vecc{g}\vecc{K}_1 \vecc{h})^T y} dy,\end{equation}
where $\vecc{\Theta}$ is, as we have remarked earlier, the time integral of the correlation function of the stationary colored noise process $\tilde{\vecc{\xi}}_t = \vecc{\sigma}_f \vecc{C}_f \vecc{\beta}_t^{(f)}$ with $\vecc{\sigma}_f$ a constant. 

We remark that if $\vecc{\Theta}_A$ is non-zero, then the heat-like functional $\mathcal{Q}^\epsilon_t$ diverges in the considered limit (since $\mathcal{Q}_t^{anom} = O(1/\epsilon^2)$ as $\epsilon \to 0$), rendering the limit non-physical. In the one-dimensional setting (where $\gamma_0 = 0$, $gh > 0$), the limit of all functionals considered is well-defined and can be expressed solely in terms of trajectory of the slow process via Stratonovich procedure. In the multi-dimensional setting, this is generally not true and, in fact, the functional might even diverge in the considered limit in the absence  of symmetry of $\vecc{K}$. In the case $\vecc{\gamma}_0 = \vecc{0}$, two sufficient condition  for $d\mathcal{W}_t'=d\mathcal{R}_t = 0$ when $\vecc{\gamma}_0 = \vecc{0}$  are:
\begin{itemize} 
\item[(1)] when the fluctuation-dissipation relation holds and the driving colored noise process is an equilibrium one (in which case $\vecc{K}_i = \vecc{C}_i \vecc{\Gamma}_i^{-1} \vecc{M}_i \vecc{C}_i^T$, $i=1,f$, is symmetric) -- this is the condition in Corollary \ref{spec1};
\item[(2)] when $\vecc{K}_1$ and $\vecc{K}_f$ are proportional to identity (but not necessarily the same), $\vecc{g}\vecc{h}$ is positive definite and commutes with $\vecc{\sigma}_f \vecc{\sigma}_f^T$. 
\end{itemize}

\subsection{The Small Mass Limit}
We introduce the scaling $m \mapsto m_0 \epsilon$ in the GLE and take the limit $\epsilon \to 0$ of the resulting equivalent rescaled SDE system:
\begin{align}
d\vecc{x}^\epsilon_t  &= \vecc{v}^\epsilon_t dt, \label{sm_gle1} \\ 
m_0 \epsilon d\vecc{v}^\epsilon_t &= \vecc{F}(t,\vecc{x}^\epsilon_t) dt - \vecc{\gamma}_0(\vecc{x}^\epsilon_t) \vecc{v}^\epsilon_t dt -  \vecc{g}(\vecc{x}^\epsilon_t) \vecc{C}_1 \vecc{y}^\epsilon_t dt +  \vecc{\sigma}_{f}(\vecc{x}^\epsilon_t) \vecc{C}_f \vecc{\beta}^{(f)\epsilon}_t dt + \vecc{\sigma}_s(\vecc{x}^\epsilon_t) \vecc{C}_s \vecc{\beta}_t^{(s)\epsilon} dt, \\
d\vecc{y}^\epsilon_t &= -\vecc{\Gamma}_1 \vecc{y}^\epsilon_t dt + \vecc{M}_1 \vecc{C}_1^T \vecc{h}(\vecc{x}^\epsilon_t) \vecc{v}^\epsilon_t dt, \\ 
d\vecc{\beta}^{(f)\epsilon}_t &= -\vecc{\Gamma}_f  \vecc{\beta}_t^{(f)\epsilon} dt + \vecc{\Sigma}_f d\vecc{W}_t^{(d_f)},  \label{sm_gle4} \\
d\vecc{\beta}^{(s)\epsilon}_t &= -\vecc{\Gamma}_s \vecc{\beta}_t^{(s)\epsilon} dt + \vecc{\Sigma}_s d\vecc{W}_t^{(d_s)}.  \label{sm_gle5}
\end{align}

The heat-like functional $\mathcal{Q}_t$ and work-like functional $\mathcal{W}_t$ satisfy the following  SDEs:
\begin{align}
d\mathcal{Q}^\epsilon_t &=  m_0 \epsilon \vecc{v}^\epsilon_t \cdot d\vecc{v}^\epsilon_t - \vecc{F}(t,\vecc{x}^\epsilon_t) \cdot d\vecc{x}^\epsilon_t , \label{sm_heat} \\ 
d\mathcal{W}^\epsilon_t &= \frac{\partial U}{\partial t} dt + \vecc{f}_{nc}(t,\vecc{x}^\epsilon_t) \cdot d\vecc{x}^\epsilon_t,  \label{sm_work}
\end{align}
where $(\vecc{x}^\epsilon_t,\vecc{v}^\epsilon_t)$ solves the SDE system \eqref{sm_gle1}-\eqref{sm_gle5}.

The dynamics in $\vecc{v}^\epsilon$ are an order of magnitude faster than those in the other variables.  Under a crucial assumption on the damping matrix $\vecc{\gamma}_0$, the limit is well-defined and we have the following results. 

\begin{cor} \label{onlysm}
Under appropriate assumptions on the initial conditions and the coefficients (i.e. Assumption \ref{ass_k}-\ref{ass_initcond}) of the pre-limit SDEs \eqref{sm_gle1}-\eqref{sm_gle5} and Assumption \ref{ass_gamma0}, the family of  processes $\vecc{x}^\epsilon_t$, satisfying the SDEs \eqref{sm_gle1}-\eqref{sm_gle5}, converges, as $\epsilon \to 0$, to the solution of the following It\^o SDE: 
\begin{align}
d\vecc{x}_t &= \vecc{\gamma}_0^{-1}(\vecc{x}_t) [\vecc{F}(t,\vecc{x}_t) -\vecc{g}(\vecc{x}_t) \vecc{C}_1 \vecc{y}_t +\vecc{\sigma}_f(\vecc{x}_t) \vecc{C}_f \vecc{\beta}_t^{(f)} + \vecc{\sigma}_s(\vecc{x}_t) \vecc{C}_s \vecc{\beta}_t^{(s)}  ] dt, \label{sm_limitx} \\
d\vecc{y}_t &= -\vecc{\Gamma}_1 \vecc{y}_t dt \nonumber \\ 
&\ \ \ \ + \vecc{M}_1 \vecc{C}_1^T \vecc{h}(\vecc{x}_t) \vecc{\gamma}_0^{-1}(\vecc{x}_t)[\vecc{F}(t,\vecc{x}_t) -\vecc{g}(\vecc{x}_t) \vecc{C}_1 \vecc{y}_t +\vecc{\sigma}_f(\vecc{x}_t) \vecc{C}_f \vecc{\beta}_t^{(f)} + \vecc{\sigma}_s(\vecc{x}_t) \vecc{C}_s \vecc{\beta}_t^{(s)}  ] dt, \\ 
d\vecc{\beta}_t^{(f)} &= -\vecc{\Gamma}_f \vecc{\beta}_t^{(f)} dt + \vecc{\Sigma}_f d\vecc{W}_t^{(d_f)}, \\
d\vecc{\beta}_t^{(s)} &= -\vecc{\Gamma}_s \vecc{\beta}_t^{(s)} dt + \vecc{\Sigma}_s d\vecc{W}_t^{(d_s)}. \label{sm_limitbeta}
\end{align}
The convergence is in the strong pathwise sense as before.
\end{cor}

Corollary \ref{onlysm} is a singular perturbation result, in a similar vein to that of \citep{hottovy2012noise,hottovy2015smoluchowski}. It effectively eliminates  description of the relatively faster velocity process from the limiting equation, keeping other processes whose time scales are slower in the equation. Because of this, the limiting equation for the position can be seen as a GLE with non-trivial non-Markovian terms. When the FDR of the second kind (i.e., Relation \ref{ass_fdr}(b)) holds, Eq. \eqref{sm_limitx}-\eqref{sm_limitbeta} simplifies to:
\begin{align}
d\vecc{x}_t &= \vecc{\gamma}_0^{-1}(\vecc{x}_t) [\vecc{F}(t,\vecc{x}_t) + \vecc{g}(\vecc{x}_t) \vecc{C}_1 \vecc{z}_t ] dt,  \\
d\vecc{z}_t &= -(\vecc{\Gamma}_1 + \vecc{M}_1 \vecc{C}_1^T \vecc{g}^T(\vecc{x}_t) \vecc{\gamma}_0^{-1}(\vecc{x}_t) \vecc{g}(\vecc{x}_t) \vecc{C}_1) \vecc{z}_t dt - \vecc{M}_1 \vecc{C}_1^T \vecc{g}^T(\vecc{x}_t) \vecc{\gamma}_0^{-1}(\vecc{x}_t)\vecc{F}(t,\vecc{x}_t) dt \nonumber \\ 
&\ \ \ \ + \vecc{\Sigma}_1 d\vecc{W}_t , 
\end{align}
where we have introduced the process $\vecc{z}_t :=  \vecc{\beta}_t - \vecc{y}_t$. The process $\vecc{z}_t$, satisfying a linear (conditional on $\vecc{x}_t$) It\^o SDE, can be thought as an effective noise term impacting the dynamics of the position process $\vecc{x}_t$ in the considered small mass limit. In contrast to the effective SDE obtained in the small mass limit of underdamped Langevin equations (i.e., those without memory and colored noise terms), the effective SDE for the position process in our GLE case is generally non-Markov.

\begin{cor}
Under the same assumptions as in Corollary \ref{onlysm}, as $\epsilon \to 0$, the family of processes $(\mathcal{W}_t^\epsilon, \mathcal{R}_t^\epsilon)$, satisfying the SDEs \eqref{sm_heat}-\eqref{sm_work}, converges to the solution of the following SDEs:
\begin{align}
d\mathcal{W}_t &= \frac{\partial U}{\partial t} dt + \vecc{f}_{nc}(t,\vecc{x}_t)  d\vecc{x}_t, \\ 
d\mathcal{R}_t &= \vecc{F}(t,\vecc{x}_t)  d\vecc{x}_t,
\end{align}
where $\vecc{x}_t$ solves the SDE \eqref{sm_limitx}. The convergence is in the strong pathwise sense as before.
\end{cor}
Note the above functionals are uniquely defined. The above corollary tells us that the limit of the dynamical process is sufficient to determine the limit of the functionals, and anomalous contributions do not appear in the small mass limit for the work and heat functionals.

\subsection{The Small Mass Limit Followed by a Markovian Limit} 
We introduce the scaling  $\vecc{\kappa}(t) \mapsto \frac{1}{\epsilon} \vecc{\kappa}\left(\frac{t}{\epsilon} \right)$ and $\vecc{R}_f(t) \mapsto \frac{1}{\epsilon} \vecc{R}_f\left(\frac{t}{\epsilon}\right) $    in the  SDEs \eqref{sm_limitx}-\eqref{sm_limitbeta}. This is the limit where a Markovian limit is taken after the small mass limit is performed on the GLE. 

The resulting rescaled SDEs for the dynamics and functionals become:
\begin{align}
d\vecc{x}^\epsilon_t &= \vecc{\gamma}_0^{-1}(\vecc{x}^\epsilon_t) [\vecc{F}(t,\vecc{x}^\epsilon_t) -\vecc{g}(\vecc{x}^\epsilon_t) \vecc{C}_1 \vecc{y}^\epsilon_t +\vecc{\sigma}_f(\vecc{x}^\epsilon_t) \vecc{C}_f \vecc{\beta}_t^{(f)\epsilon} + \vecc{\sigma}_s(\vecc{x}_t^\epsilon) \vecc{C}_s \vecc{\beta}_t^{(s)\epsilon}  ] dt, \label{sm_prelimitx} \\
\epsilon d\vecc{y}^\epsilon_t &= -\vecc{\gamma}_1(\vecc{x}^\epsilon_t) \vecc{y}^\epsilon_t dt \nonumber \\ 
&\ \  \ \ + \vecc{M}_1 \vecc{C}_1^T \vecc{h}(\vecc{x}^\epsilon_t) \vecc{\gamma}_0^{-1}(\vecc{x}^\epsilon_t)[\vecc{F}(t,\vecc{x}^\epsilon_t) +\vecc{\sigma}_f(\vecc{x}^\epsilon_t) \vecc{C}_f \vecc{\beta}_t^{(f)\epsilon} + \vecc{\sigma}_s(\vecc{x}^\epsilon_t) \vecc{C}_s \vecc{\beta}_t^{(s)\epsilon}  ] dt, \\ 
\epsilon d\vecc{\beta}_t^{(f)\epsilon} &= -\vecc{\Gamma}_f \vecc{\beta}_t^{(f)\epsilon} dt + \vecc{\Sigma}_f d\vecc{W}_t^{(d_f)}, \\
d\vecc{\beta}_t^{(s)\epsilon} &= -\vecc{\Gamma}_s \vecc{\beta}_t^{(s)\epsilon} dt + \vecc{\Sigma}_s d\vecc{W}_t^{(d_s)}, \label{sm_prelimitbeta} \\ 
d\mathcal{W}^\epsilon_t &= \frac{\partial U}{\partial t} dt + \vecc{f}_{nc}(t,\vecc{x}^\epsilon_t)  d\vecc{x}^\epsilon_t, \label{smm_work} \\ 
d\mathcal{R}^\epsilon_t &= \vecc{F}(t,\vecc{x}^\epsilon_t) d\vecc{x}^\epsilon_t, \label{smm_heat}
\end{align}
where $\vecc{\gamma}_1 = \vecc{\Gamma}_1 + \vecc{M}_1 \vecc{C}_1^T \vecc{h} \vecc{\gamma}_0^{-1} \vecc{g} \vecc{C}_1$. 

\begin{cor} \label{smm}
Under appropriate assumptions on the initial conditions and the coefficients (i.e. Assumption \ref{ass_k}-\ref{ass_initcond}) of the pre-limit SDEs \eqref{sm_prelimitx}-\eqref{sm_prelimitbeta} and Assumption \ref{ass_gamma0}, the family of  processes $\vecc{x}^\epsilon_t$, satisfying the SDEs \eqref{sm_prelimitx}-\eqref{sm_prelimitbeta}, converges, as $\epsilon \to 0$, to the solution of the following It\^o SDE: 
\begin{align}
d\vecc{x}_t &= \vecc{\gamma}_2^{-1}(\vecc{x}_t) [ \vecc{F}(t,\vecc{x}_t)  + \vecc{\sigma}_s(\vecc{x}_t) \vecc{C}_s \vecc{\beta}_t^{(s)}  ] dt + \vecc{\gamma}_2^{-1}(\vecc{x}_t) \vecc{\sigma}_f \vecc{C}_f \vecc{\Gamma}_f^{-1} \vecc{\Sigma}_f  d\vecc{W}_t^{(d_f)} + \vecc{S}(\vecc{x}_t) dt, \label{smm_limitx} \\
d\vecc{\beta}_t^{(s)} &= -\vecc{\Gamma}_s \vecc{\beta}_t^{(s)} dt + \vecc{\Sigma}_s d\vecc{W}_t^{(d_s)}, \label{smm_limitbeta}
\end{align}
where $\vecc{\gamma}_2^{-1} = \vecc{\gamma}_0^{-1}(\vecc{I} - \vecc{g} \vecc{C}_1 \vecc{\gamma}_1^{-1} \vecc{M}_1 \vecc{C}_1^T \vecc{h} \vecc{\gamma}_0^{-1}) $, $\vecc{\gamma}_1 = \vecc{\Gamma}_1 + \vecc{M}_1 \vecc{C}_1^T \vecc{h} \vecc{\gamma}_0^{-1} \vecc{g} \vecc{C}_1$, and 
\begin{align}
S^i &= \frac{\partial R^{ij}}{\partial x^l} T^{jl}.
\end{align} 
In the above 
\begin{equation}
\vecc{R} = -\vecc{\gamma}_0^{-1}[\vecc{g}\vecc{C}_1\vecc{\gamma}_1^{-1} \hspace{0.8cm} \vecc{g}\vecc{C}_1\vecc{\gamma}_1^{-1}(\vecc{M}_1\vecc{C}_1^T\vecc{h}\vecc{\gamma}_0^{-1}\vecc{\sigma}_f\vecc{C}_f)\vecc{\Gamma}_f^{-1} - \vecc{\sigma}_f\vecc{C}_f\vecc{\Gamma}_f^{-1}], 
\end{equation}
\begin{equation}
\vecc{T} = (-\vecc{J}_{11}\vecc{C}_1^T \vecc{g}^T \vecc{\gamma}_0^{-T} + \vecc{J}_{12}\vecc{C}_f^T \vecc{\sigma}_f^T \vecc{\gamma}_0^{-T}, -\vecc{J}_{12}^T\vecc{C}_1^T \vecc{g}^T \vecc{\gamma}_0^{-T}+\vecc{M}_f \vecc{C}_f^T \vecc{\sigma}_f^T \vecc{\gamma}_0^{-T}),
\end{equation}
where $\vecc{J}_{11}$ and $\vecc{J}_{12}$ solve the matrix equations:
\begin{align}
\vecc{\gamma}_1 \vecc{J}_{12} + \vecc{J}_{12} \vecc{\Gamma}_f^T &= \vecc{M}_1 \vecc{C}_1^T \vecc{h} \vecc{\gamma}_0^{-1} \vecc{\sigma}_f \vecc{C}_f \vecc{M}_f, \label{meq1} \\ 
\vecc{\gamma}_1 \vecc{J}_{11} + \vecc{J}_{11} \vecc{\gamma}_1^T &= \vecc{M}_1 \vecc{C}_1^T \vecc{h} \vecc{\gamma}_0^{-1} \vecc{\sigma}_f \vecc{C}_f \vecc{J}_{12}^T + \vecc{J}_{12} ( \vecc{M}_1 \vecc{C}_1^T \vecc{h} \vecc{\gamma}_0^{-1} \vecc{\sigma}_f \vecc{C}_f)^T.  \label{meq2}
\end{align}
The convergence is in the strong pathwise sense as before.
\end{cor}

Similar to Corollary \ref{follow_limit_res1}, Corollary \ref{smm} says that taking the small mass limit and then the Markovian limit on the GLE \eqref{genle_general} leads to  effective  dynamics for the position that can be described by a Langevin equation driven by a white noise and a colored noise $\vecc{\beta}_t^{(s)}$, but with an additional drift correction term $\vecc{S}$ and without a memory term.  In other words, the effective dynamics can be equivalently described by:
\begin{align}
\dot{\vecc{x}}_t &= \vecc{\gamma}_2^{-1}(\vecc{x}_t) (\vecc{F}(t,\vecc{x}_t)  + \vecc{\sigma}_f \vecc{C}_f \vecc{\Gamma}_f^{-1} \vecc{\Sigma}_f \dot{\vecc{W}}^{(d_f)}_t + \vecc{\sigma}_s(\vecc{x}_t) \vecc{C}_s \vecc{\beta}_t^{(s)} ) + \vecc{S}(\vecc{x}_t).
\end{align}
However, comparing the above equation to Eq. \eqref{88}, we see that the procedures of taking the small mass limit and taking the Markovian limit do not commute, regardless of whether the FDR holds or not, and hence the order of taking limits matters here.

\begin{cor}
Under the same assumptions as in Corollary  \ref{smm}, as $\epsilon \to 0$,  the family of processes $(\mathcal{W}_t^\epsilon, \mathcal{R}_t^\epsilon)$, satisfying the SDEs \eqref{smm_heat}-\eqref{smm_work}, converges to the solution of the following SDEs:
\begin{align}
d\mathcal{W}_t &= \frac{\partial U}{\partial t} dt + \vecc{f}_{nc}(t,\vecc{x}_t) \circ  d\vecc{x}_t + d\mathcal{W}_t', \\ 
d\mathcal{R}_t &= \vecc{F}(t,\vecc{x}_t) \circ  d\vecc{x}_t + d\mathcal{R}_t', \\
d\mathcal{W}_t' &= [\vecc{\nabla} \cdot (\vecc{f}_{nc}^T(t,\vecc{x}_t) \vecc{\Phi}(\vecc{x}_t) \vecc{\mu}_{A}^T(\vecc{x}_t) \vecc{\Phi}^T(\vecc{x}_t)) - \vecc{f}_{nc}^T(t,\vecc{x}_t) \vecc{\nabla} \cdot (\vecc{\Phi}(\vecc{x}_t) \vecc{\mu}_A^T(\vecc{x}_t) \vecc{\Phi}^T(\vecc{x}_t) ) ] dt, \\
d\mathcal{R}_t'&=  [\vecc{\nabla} \cdot (\vecc{F}^T(t,\vecc{x}_t) \vecc{\Phi}(\vecc{x}_t) \vecc{\mu}_{A}^T(\vecc{x}_t) \vecc{\Phi}^T(\vecc{x}_t)) - \vecc{F}^T(t, \vecc{x}_t) \vecc{\nabla} \cdot (\vecc{\Phi}(\vecc{x}_t) \vecc{\mu}_A^T(\vecc{x}_t) \vecc{\Phi}^T(\vecc{x}_t) ) ] dt,
\end{align}
where 
$\vecc{\Phi} = \vecc{\gamma}_0^{-1}[-\vecc{g}\vecc{C}_1 \ \ \ \vecc{\sigma}_f \vecc{C}_f]$, $\vecc{\mu}_A$ is the antisymmetric part of the matrix 
\begin{equation}
\vecc{\mu} = 
\begin{bmatrix}
\vecc{\gamma}_1^{-1}(\vecc{J}_{11} + \vecc{M}_1 \vecc{C}_1^T \vecc{h} \vecc{\gamma}_0^{-1} \vecc{\sigma}_f \vecc{C}_f \vecc{\Gamma}_f^{-1} \vecc{J}_{12}^T) & \vecc{\gamma}_1^{-1}(\vecc{J}_{12} + \vecc{M}_1 \vecc{C}_1^T \vecc{h} \vecc{\gamma}_0^{-1} \vecc{\sigma}_f \vecc{C}_f \vecc{\Gamma}_f^{-1} \vecc{M}_f) \\
\vecc{\Gamma}_f^{-1} \vecc{J}_{12}^T & \vecc{\Gamma}_f^{-1} \vecc{M}_f \\ 
\end{bmatrix},
\end{equation}
with $\vecc{J}_{11}$ and $\vecc{J}_{12}$ satisfying \eqref{meq1}-\eqref{meq2}, and $\vecc{x}_t$ solves the SDE \eqref{smm_limitx}. The convergence is in the strong pathwise sense as before.
\end{cor}

Comparing the above result to that in Corollary \ref{follow_limit_res}, we see that the procedures of taking the small mass limit and taking the Markovian limit again do not commute, regardless of whether the FDR holds or not and whether $\vecc{\Theta}$ is symmetric or not, and hence the order of taking limits also matters for the functionals.

\subsection{A Joint Markovian and Small Mass Limit} \label{sect_limit3} 
We introduce the scaling  $\vecc{\kappa}(t) \mapsto \frac{1}{\epsilon} \vecc{\kappa}\left(\frac{t}{\epsilon} \right)$ and $\vecc{R}_f(t) \mapsto \frac{1}{\epsilon} \vecc{R}_f\left(\frac{t}{\epsilon}\right) $, $m \mapsto m_0 \epsilon$ in the GLE \eqref{genle_general}. This is the limit where the inertial time scale, the memory time scale and some noise correlation time scales of the system tend to zero at the same rate. This will provide a further coarse-grained model compared to the Markovian limit and therefore more information will be lost in the limit.   We remark that the small mass limit of our GLE  is generally not well-defined (unless $\vecc{\gamma}_0 > \vecc{0}$) and leads to the interesting phenomenon of anomalous gap of the particle's mean-squared displacement  \citep{mckinley2009transient,indei2012treating,cordoba2012elimination}. 

Introducing the auxiliary variable $\vecc{y}_t$ as before, the resulting rescaled GLE can then be studied as the following SDE system for the Markov process $(\vecc{x}^\epsilon_t, \vecc{v}^\epsilon_t, \vecc{y}^\epsilon_t, \vecc{\beta}^{(f)\epsilon}_t, \vecc{\beta}_t^{(s)\epsilon})$: 
\begin{align}
d\vecc{x}^\epsilon_t  &= \vecc{v}^\epsilon_t dt, \label{joint_gle1} \\ 
m_0 \epsilon d\vecc{v}^\epsilon_t &= \vecc{F}(t,\vecc{x}^\epsilon_t) dt -  \vecc{\gamma}_0(\vecc{x}^\epsilon_t) \vecc{v}^\epsilon_t dt - \vecc{g}(\vecc{x}^\epsilon_t) \vecc{C}_1  \vecc{y}^\epsilon_t dt +  \vecc{\sigma}_f(\vecc{x}^\epsilon_t) \vecc{C}_f \vecc{\beta}_t^{(f)\epsilon} dt + \vecc{\sigma}_s(\vecc{x}^\epsilon_t) \vecc{C}_s \vecc{\beta}_t^{(s)\epsilon} dt , \\
\epsilon d\vecc{y}^\epsilon_t &= -\vecc{\Gamma}_1 \vecc{y}^\epsilon_t dt + \vecc{M}_1 \vecc{C}_1^T \vecc{h}(\vecc{x}^\epsilon_t) \vecc{v}^\epsilon_t dt, \\ 
\epsilon d\vecc{\beta}^{(f)\epsilon}_t &= -\vecc{\Gamma}_f  \vecc{\beta}^{(f)\epsilon}_t dt + \vecc{\Sigma}_f d\vecc{W}_t^{(d_f)},  \label{joint_gle4} \\
d\vecc{\beta}^{(s)\epsilon}_t &= -\vecc{\Gamma}_s \vecc{\beta}^{(s)\epsilon}_t dt + \vecc{\Sigma}_s d\vecc{W}_t^{(d_s)}. \label{joint_gle5}
\end{align}

The heat $\mathcal{Q}^\epsilon_t$ and work $\mathcal{W}^\epsilon_t$ satisfy the following SDEs:
\begin{align}
d\mathcal{Q}^\epsilon_t &=  m_0  \epsilon \vecc{v}^\epsilon_t \cdot d\vecc{v}_t^\epsilon - \vecc{F}(t,\vecc{x}^\epsilon_t) \cdot d\vecc{x}^\epsilon_t, \\ 
d\mathcal{W}^\epsilon_t &= \frac{\partial U}{\partial t} dt + \vecc{f}_{nc}(t,\vecc{x}^\epsilon_t) \cdot d\vecc{x}^\epsilon_t, 
\end{align}
where $(\vecc{x}^\epsilon_t, \vecc{v}^\epsilon_t)$ solves the SDE system \eqref{joint_gle1}-\eqref{joint_gle5}.

The dynamics in $\vecc{v}^\epsilon$, $\vecc{y}^\epsilon$ and $\vecc{\beta}^{(f)\epsilon}$ are an order of magnitude faster than those in $\vecc{x}^\epsilon$, $\vecc{\beta}^{(s)\epsilon}$, $\mathcal{Q}^\epsilon$ and $\mathcal{W}^\epsilon$. 

Consider the following system of five matrix equations for $\vecc{J}_{11} = \vecc{J}_{11}^T$, $\vecc{J}_{21}=\vecc{J}_{12}^T $ and $\vecc{J}_{31} = \vecc{J}_{13}^T$ (c.f. \citep{Lim2018}):
\begin{align}
\vecc{\gamma}_0 \vecc{J}_{11} + \vecc{J}_{11} \vecc{\gamma}_0^T +  \vecc{g} \vecc{C}_1 \vecc{J}_{12}^T + \vecc{J}_{12} \vecc{C}_1^T \vecc{g}^T &= \vecc{\sigma}_f \vecc{C}_f \vecc{J}_{13}^T + \vecc{J}_{13} \vecc{C}_f^T \vecc{\sigma}_f^T, \label{gen_system} \\
m_0 \vecc{J}_{11} \vecc{h}^T \vecc{C}_1 \vecc{M}_1 + \vecc{\sigma}_f \vecc{C}_f \vecc{J}_{23}^T &=  \vecc{g} \vecc{C}_1 \vecc{J}_{22} + m_0 \vecc{J}_{12} \vecc{\Gamma}_1^T + \vecc{\gamma}_0 \vecc{J}_{12},\\
\vecc{\gamma}_0 \vecc{J}_{13} +  \vecc{g} \vecc{C}_1 \vecc{J}_{23} + m_0 \vecc{J}_{13} \vecc{\Gamma}_f^T &= \vecc{\sigma}_f \vecc{C}_f \vecc{M}_f,  \\ 
\vecc{M}_1 \vecc{C}_1^T \vecc{h} \vecc{J}_{12} + \vecc{J}_{12}^T \vecc{h}^T \vecc{C}_1 \vecc{M}_1 &= \vecc{\Gamma}_1 \vecc{J}_{22} + \vecc{J}_{22} \vecc{\Gamma}_1^T, \\ 
\vecc{M}_1 \vecc{C}_1^T \vecc{h} \vecc{J}_{13} &=  \vecc{\Gamma}_1 \vecc{J}_{23} + \vecc{J}_{23} \vecc{\Gamma}_f^T.  \label{gen_end} 
\end{align}

We write  $\mathcal{Q}^\epsilon_t =  \frac{m_0}{2}  \epsilon |\vecc{v}^\epsilon_t|^2 - \frac{m_0}{2}  \epsilon |\vecc{v}^\epsilon_0|^2  -\mathcal{R}^\epsilon_t$, where $\mathcal{R}^\epsilon_t = \int_0^t \vecc{F}(s,\vecc{x}^\epsilon_s) \cdot d\vecc{x}^\epsilon_s$. We expect that as $\epsilon \to 0$,  the kinetic energy terms are of $O(1)$ and they tend to $\overline{\frac{m_0}{2}  |\vecc{v}_t|^2 - \frac{m_0}{2}  |\vecc{v}_0|^2}$, where the overline denotes average with respect to the invariant density of the stationary fast process (at a given slow ones), which is mean zero Gaussian with covariance matrix $\vecc{J}_{11}$.  Therefore, to study the asymptotic behavior of $\mathcal{Q}^\epsilon_t$ in the considered limit, it suffices to investigate the asymptotic behavior of $\mathcal{R}^\epsilon_t$. 

One then has the following results.

\begin{cor} \label{jointlimit_res1}
The family of processes $\vecc{x}^\epsilon_t$, satisfying the SDEs \eqref{joint_gle1}-\eqref{joint_gle5}, converges, as $\epsilon \to 0$, to the solution of the following It\^o SDE: 
\begin{align}
d\vecc{x}_t &= \vecc{\Gamma}^{-1}(\vecc{x}_t) (\vecc{F}(t,\vecc{x}_t) + \vecc{\sigma}_s(\vecc{x}_t) \vecc{C}_s \vecc{\beta}_t^{(s)} ) dt +  \vecc{S}(\vecc{x}_t) dt + \vecc{\Gamma}^{-1}(\vecc{x}_t)\vecc{\Sigma}(\vecc{x}_t) d\vecc{W}_t^{(d_f)}, \label{jointlimit_eqx} \\
d\vecc{\beta}_t^{(s)} &= -\vecc{\Gamma}_s \vecc{\beta}_t^{(s)} dt + \vecc{\Sigma}_s d\vecc{W}_t^{(d_s)}, \label{jointlimit_eqbeta}
\end{align}
where $\vecc{\Gamma} = \vecc{\gamma}_0 + \vecc{g} \vecc{K}_1 \vecc{h}$, $\vecc{\Sigma} = \vecc{\sigma}_f \vecc{C}_f \vecc{\Gamma}_f^{-1} \vecc{\Sigma}_f$, and $\vecc{S}$ is the noise-induced drift whose expression is given by:
\begin{align}
\vecc{S} &= \vecc{\nabla} \cdot \left(\vecc{\Gamma}^{-1}(m_0 \vecc{J}_{11} - \vecc{g} (\vecc{C}_1 \vecc{\Gamma}_1^{-1} \vecc{J}_{21})^T + \vecc{\sigma}_f(\vecc{C}_f \vecc{\Gamma}_f^{-1} \vecc{J}_{31})^T ) \right) \nonumber \\ 
&\ \ \ + \vecc{\Gamma}^{-1}\left( \vecc{g} \vecc{\nabla}\cdot ((\vecc{C}_1 \vecc{\Gamma}_1^{-1} \vecc{J}_{21})^T) - \vecc{\sigma}_f \vecc{\nabla} \cdot ((\vecc{C}_f \vecc{\Gamma}_f^{-1} \vecc{J}_{31})^T) - m_0 \vecc{\nabla} \cdot \vecc{J}_{11}\right),  \label{nid1}
\end{align}
 where the $\vecc{J}_{ij}$ solve the system of matrix equations \eqref{gen_system}-\eqref{gen_end}. The convergence is in the strong pathwise sense as before.
\end{cor}

The results in Corollary \ref{jointlimit_res1} have the same structure as those in Corollary \ref{follow_limit_res1}, but with different drift correction term (see also the discussion following Corollary \ref{follow_limit_res1}). The presence of the noise-induced drift $\vecc{S}$, due to the state-dependence of the coefficients $\vecc{g}$, $\vecc{h}$ and $\vecc{\sigma}_f$, and the dependence of its formula on the limiting procedure taken implies that the elimination of the fast degrees of freedom needs to be done carefully and naive procedure could lead to inconsistent result. In the special case when the FDR relation of the second kind holds, the noise-induced drift in Eq. \eqref{nid1} simplifies to $\vecc{S} = m_0 (\vecc{\nabla} \cdot (\vecc{\Gamma}^{-1} \vecc{J}_{11}) - \vecc{\Gamma}^{-1} \vecc{\nabla} \cdot \vecc{J}_{11})$ (see also Corollary 2 in \citep{Lim2018}).

\begin{cor} \label{joint_limit_res}
 Let $\vecc{\lambda}_A$ denote the antisymmetric part of $\vecc{\lambda} = -m_0 \vecc{\Gamma}^{-1} \vecc{J}_{11}+ \vecc{\Gamma}^{-1} \vecc{g} \vecc{C}_1 \vecc{\Gamma}_1^{-1} \vecc{J}_{21} -  \vecc{\Gamma}^{-1} \vecc{\sigma}_f \vecc{C}_f \vecc{\Gamma}_f^{-1} \vecc{J}_{31}$,  where the $\vecc{J}_{ij}$ solve the system of matrix equation \eqref{gen_system}-\eqref{gen_end}. 
  
The family of processes $(\mathcal{W}^\epsilon_t, \mathcal{R}^\epsilon_t)$ converges, as $\epsilon \to 0$, to the solution of the following SDEs: 
\begin{align}
d\mathcal{W}_t &= \frac{\partial U}{\partial t} dt + \vecc{f}_{nc}(t,\vecc{x}_t) \circ  d\vecc{x}_t + dW_t^{anom}, \label{jointlimit_eqW} \\
d\mathcal{R}_t &= \vecc{F}(t, \vecc{x}_t) \circ d\vecc{x}_t + d\mathcal{R}_t^{anom},
\end{align}
where 
\begin{align}
d\mathcal{W}_t^{anom} &= [\vecc{\nabla} \cdot (\vecc{f}_{nc}^{T}(t,\vecc{x}_t) \vecc{\lambda}_A(\vecc{x}_t)) -  \vecc{f}_{nc}^{T}(t,\vecc{x}_t)  \vecc{\nabla} \cdot \vecc{\lambda}_A(\vecc{x}_t)] dt, \\
d\mathcal{R}_t^{anom} &= [\vecc{\nabla} \cdot (\vecc{F}^{T}(t,\vecc{x}_t) \vecc{\lambda}_A(\vecc{x}_t)) -  \vecc{F}^{T}(t,\vecc{x}_t)  \vecc{\nabla} \cdot \vecc{\lambda}_A(\vecc{x}_t)] dt,
\end{align}
and $\vecc{x}_t$ solves the SDE \eqref{jointlimit_eqx}-\eqref{jointlimit_eqbeta}. The convergence is in the strong pathwise sense as before.  
\end{cor}

\begin{cor}
$d\mathcal{W}_t^{anom}=d\mathcal{R}_t^{anom}=0$ when one of the following conditions holds:
\begin{itemize}
\item[(i)] $\vecc{\gamma}_0$, $\vecc{g}$, $\vecc{h}$, $\vecc{\sigma}_f$, $\vecc{C}_i$, $\vecc{M}_i$, $\vecc{\Gamma}_i$ ($i=1,f$) are diagonal; 
\item[(ii)] $\vecc{\gamma}_0 = \vecc{0}$, the fluctuation-dissipation relation of the second kind holds, and $\vecc{\Gamma}^{-1} \vecc{\sigma}_f \vecc{K}_f^T \vecc{\sigma}_f^T$ is symmetric. 
\end{itemize}  
\end{cor}

Note that the FDR relation of the second kind alone is not enough to guarantee $d\mathcal{W}_t^{anom}=d\mathcal{R}_t^{anom}=0$. 

In contrast to the Markovian limit case, it is generally not possible to express both the work and heat functional in terms of trajectory of the effective slow process without additional drift terms. This is possible for the work functional in the case where $\vecc{f}_{nc}$ is independent of position. Also, the matrix $\vecc{\lambda}$  loses the meaning as the time integral of the correlation function of a physical noise process. 

We next discuss the above results in the case $\vecc{\gamma}_0 = \vecc{0}$. 
The limiting expression for $\mathcal{W}_t$ and $\mathcal{R}_t$ can be expressed in terms of trajectory of the slow process via Stratonovich discretization if and only if $\vecc{\lambda}_A$ vanishes.  
In the one-dimensional setting, the Stratonovich procedure is justified even if the fluctuation-dissipation relation is violated. However, in contrast to the results obtained for the Markovian limit, a stricter condition is needed for $\vecc{\lambda}_A$ to vanish in the general multi-dimensional case.  Whether $\vecc{\lambda}_A$ vanishes or not  is not entirely attributed to the symmetry associated with the noise term, but it also depends on the properties of the memory function as well as the coefficients $\vecc{g}$, $\vecc{h}$ and $\vecc{\sigma}_f$.  The unifying message in the above discussion is that, in the multidimensional setting, higher level of coarse-graining or model reduction often leads to justification of use of Stratonovich procedure in defining thermodynamic functionals using equations for the effective dynamics for a smaller, more restricted class of systems. In the special one-dimensional setting, the Stratonovich procedure is always justified.

\subsection{Further Discussions}

%different levels of symmetries and fluctuation-dissipation relation
%commutativity -- state-dependent coeff

We have considered the joint Markovian and small mass limit of the GLE (Procedure (5.5)) in the previous subsection, as well as the procedure where the small mass limit is taken after the Markovian limit is taken here (Procedure (5.2)).  A natural question is how do the effective equations obtained via these two limiting procedures compare. 

To allow the comparison, we assume that $\vecc{\gamma}_0 = \vecc{0}$ and $\vecc{\Theta}_A = \vecc{0}$ in the following discussion.  First, note that the solution of \eqref{jointlimit_eqx} coincides, in law, with that of \eqref{difflimit_eqx} if and only if the noise-induced drifts $\vecc{S}$ (in \eqref{nid1})  and $\vecc{H}$ (in \eqref{nid2}) coincide.  Second, the work functionals, satisfying \eqref{jointlimit_eqW} and \eqref{difflimit_eqW}  respectively, coincide, if in addition, $\mathcal{W}_t^{anom} = \mathcal{W}_t'$, i.e. if and only if  $\vecc{\lambda}_A = \vecc{K}_A^T$. These occur, for instance, in the very special case of one dimensions where the fluctuation-dissipation relation of the second kind, i.e. $\vecc{g} = \vecc{h}^T = \vecc{\sigma}_f$ and $\vecc{R}_f(t) = \vecc{\kappa}(t)$ holds. In the general case, these need not hold, even if the fluctuation-dissipation relation of the second kind holds. Indeed, in this case (i.e., when the FDR of the second kind holds)  Eq. \eqref{gen_system} reduces to 
\begin{equation}
\vecc{g} \vecc{C}_1(\vecc{J}_{12}-\vecc{J}_{13})^T + (\vecc{J}_{12}-\vecc{J}_{13})(\vecc{g}\vecc{C}_1)^T = \vecc{0},
\end{equation}
and so $\vecc{J}_{12} = \vecc{J}_{13}$ and $ \vecc{J}_{21} = \vecc{J}_{31}$. Therefore, the noise-induced drift in Eq. \eqref{nid1} simplifies to $\vecc{S} = m_0 (\vecc{\nabla} \cdot (\vecc{\Gamma}^{-1} \vecc{J}_{11}) - \vecc{\Gamma}^{-1} \vecc{\nabla} \cdot \vecc{J}_{11})$. However, $m_0 \vecc{J}_{11}$ needs not equal the $\vecc{J}$ defined in Eq. \eqref{nid2}, unless we are in the special case  where all the variables commute (for instance, in one dimensions), in which case we can solve for $\vecc{J}$ explicitly and obtain $m_0 \vecc{J}_{11} = \vecc{J} = \vecc{I}$. This highlights the fact that these two procedures of taking limit do not generally give us the same effective dynamics. The key message here is that the manner in which the limits are taken matters and hold important consequences for the considered physical systems, so one should be careful when attempting to take various limits in stochastic thermodynamics.

Similar, albeit slightly more tedious, comparison can also be performed for the results obtained via Procedure (5.5) and those via Procedure (5.4).  In general, convergence of the dynamical and functional paths depends on regularity of the approximating sequence. Different homogenization procedures give rise to approximating sequences of different regularity and thus different limiting behavior where different forms of area anomaly appear, so the commutativity of the procedures is not guaranteed unless one restricts to special cases -- these cases invoke symmetry of the Onsager matrix associated with the fast dynamics as well as the relation between dissipation and fluctuation driving the fast dynamics. 

%special/sufficient conditions
%thermophoresis when T dep on x
%figure showing relationship between results

\section{Conclusions} \label{sect_concl}
We have explored and performed various multiple time scale  analysis (homogenization) for a class of generalized Langevin dynamics together with the stochastic processes describing the heat-like and work-like functionals in stochastic thermodynamics. We have addressed and discussed the important problem of justifying the use of Stratonovich convention in the definition of these functionals in the situations where there exists wide separation of time scales of various levels in the systems. We find that, unless certain symmetry is present in the GLE system, it is generally not possible to express the effective evolution of these functionals solely in terms of trajectory of the effective process describing the system dynamics via the standard Stratonovich convention, and additional information of the full process is needed to do so. 

Depending on the level of coarse graining, one needs to impose appropriate symmetry conditions in such a way that the area anomaly, encoded by the antisymmetric part of the Onsager matrix associated with the fast dynamics, vanishes, in order to make this possible. In the  case where these functionals are thermodynamic, the absence of these symmetry conditions  gives rise to anomalous thermodynamics in the homogenized systems. Our results can be applied to concrete physical systems, including the ones described in Appendix \ref{sect_eg}, in various time scale separation scenarios.   \\

%Lastly, we discuss two interesting future directions. So far we have been focusing on finite-dimensional GLE systems. It would be interesting to study homogenization for infinite-dimensional GLE systems, in which case the memory function decays as a power law, and the functionals along the trajectories of these systems. We have also been restricting our study to the case of infinite time scale separation limit.  Finite time scale limit -> corrections in term of higher order correlation functions -> rough paths  \\

%Ideas for extensions (but will not be pursued):
%\begin{itemize}
%/item Anomalous diffusion -- partial Markovian limit %(superdiffusion), Markovian limit (subdiffusion) 
%\item Only small mass limit 
%\item connection to renormalized SDEs and rough paths
%\item deterministic approximation of stochastic dynamics and functionals - presence of anomaly even for state-independent coefficients case
%\end{itemize}

\noindent {\bf Acknowledgements.}
The author gratefully acknowledges the support provided by the NORDITA Fellowship during the year 2018-2021. He is also grateful to  Peter H{\"a}nggi, Ralf Eichhorn and Stefano Bo for critical reading of the manuscript and suggestions for improving an earlier version of the manuscript. \\

%---------------------------------------------
%\section*{Appendices}
\newpage
\appendix

\section{Examples of GLE Systems in Nonequilibrium Statistical Mechanics} \label{sect_eg}
In this section, we provide three examples of physical system that can be modeled by (special cases of) the GLEs studied in this paper.

\ex {\bf A Brownian particle in a temperature gradient.}
We consider a Brownian particle immersed in a nonequilibrium heat bath where a temperature gradient is present.  For this system, the temperature of the heat bath varies with the position of the particle and a generalized fluctuation-dissipation relation holds. We model the system by the GLE defined in Section \ref{sect_gles}, with $\vecc{\gamma}_0 = \vecc{0}$, $\vecc{\sigma}_0 = \vecc{0}$,   $\vecc{g}(\vecc{x})=\vecc{h}(\vecc{x})=\sqrt{\gamma(\vecc{x})}\vecc{I}$ and $\vecc{\sigma}(\vecc{x}) = \sqrt{k_B T(\vecc{x}) \gamma(\vecc{x})} \vecc{I}$, where $\vecc{x} \in \RR^d$ ($d=1,2,3$),  $\gamma>0$ is a scalar function, $k_B$ is the Boltzmann constant, $T$ is the state-dependent temperature  of the bath and  $\vecc{I}$ is the identity matrix. The resulting GLE is then: 
%\begin{align}
%d\vecc{x}_t &= \vecc{v}_t dt, \\
%m d\vecc{v}_t &= \vecc{F}(t, \vecc{x}_t) dt -\vecc{\sigma}%(\vecc{x}_t) \left(\int_0^t \kappa(t-s) \vecc{\sigma}%^T(\vecc{x}_s) \vecc{v}_s ds\right) dt + \sqrt{2 k_B %T(\vecc{x}_t)} \vecc{\sigma}(\vecc{x}_t) \vecc{\xi}_t dt,
%\end{align}
%The second model is a special case of the first model, where %$\vecc{\sigma}(\vecc{x}_t) = \sqrt{\gamma(\vecc{x}_t)} %\vecc{I}$, and it is described by:
\begin{align}
d\vecc{x}_t &= \vecc{v}_t dt, \label{eg11} \\
m d\vecc{v}_t &= \vecc{F}(t, \vecc{x}_t) dt -\sqrt{\gamma(\vecc{x}_t)} \left(\int_0^t  \vecc{\kappa}(t-s) \sqrt{\gamma(\vecc{x}_s)} \vecc{v}_s ds\right) dt + \sqrt{ k_B T(\vecc{x}_t) \gamma(\vecc{x}_t)} \vecc{\xi}_t dt, \label{eg12}
\end{align}
where $\vecc{\xi}_t \in \RR^d$ is a mean-zero, stationary Gaussian colored noise with covariance function equals to $\vecc{\kappa}(t)$.  The memory function and the colored noise are defined in a similar way as before: $\vecc{\kappa}(t)  = \vecc{C}_1 e^{-\vecc{\Gamma}_1 t} \vecc{M}_1 \vecc{C}_1^T$  and $\vecc{\xi}_t = \vecc{C}_1 \vecc{\beta}_t$, where $d\vecc{\beta}_t = -\vecc{\Gamma}_1 \vecc{\beta}_t dt + \vecc{\Sigma}_1 d\vecc{W}_t$. The above model has been used to study the phenomena of thermophoresis in \citep{Lim2018} (see also the discussions and references related to the GLE  \eqref{eg11}-\eqref{eg12} there).

\ex {\bf Active matter systems with spatially inhomogeneous activity.}
We consider a small system in an equilibrium (passive) heat bath at the constant temperature $T$ subject to an external force field described by $\vecc{F}(t,\vecc{x}) = -\vecc{\nabla}_{\vecc{x}}U(t,\vecc{x}) + \vecc{f}_{nc}(t,\vecc{x})$ and an active force field described by $\vecc{\sigma}_a(\vecc{x}) \vecc{\eta}$, where $\vecc{x} \in \RR^d$ ($d=1,2,3$), $U$ is the potential, $\vecc{f}_{nc}$ is a non-conservative force field, $\vecc{\sigma}_a: \RR^d \to \RR^{d \times a}$ is a state-dependent coefficient, and $\vecc{\eta} \in \RR^a$ is a mean-zero stationary Ornstein-Uhlenbeck process. We model this system by the GLE in Section \ref{sect_gles} with $\vecc{\gamma}_0 = \vecc{0}$, $\vecc{\sigma}_0 = \vecc{0}$,  $\vecc{g} = \vecc{h}^T = \vecc{\sigma}_p \in \RR^{d \times d_1}$ (constant matrix), $\vecc{\sigma}(\vecc{x}) =  [\sqrt{k_B T} \vecc{\sigma}_p \ \ \vecc{\sigma}_a(\vecc{x})] \in \RR^{d \times (d_1 + a)}$, $\vecc{\xi}_t =\vecc{C}_2 \vecc{\beta}_t$, with $\vecc{C}_2 = \vecc{I}$, $\vecc{\beta}_t = (\vecc{\zeta}_t, \vecc{\eta}_t) \in \RR^{d_1 + a}$, $\vecc{\zeta}_t  = \vecc{C}_p \vecc{\theta_t}$. More precisely:
\begin{align}
d\vecc{x}_t &= \vecc{v}_t dt, \\
m d\vecc{v}_t &= \vecc{F}(t,\vecc{x}_t) dt - \vecc{\sigma}_p \left(\int_0^t \vecc{\kappa}(t-s) \vecc{\sigma}_p^T \vecc{v}_s ds\right) dt + \sqrt{k_B T} \vecc{\sigma}_p \vecc{C}_p \vecc{\theta}_t dt +   \vecc{\sigma}_a(\vecc{x}_t) \vecc{\eta}_t dt, \\ 
d\vecc{\theta}_t &= -\vecc{\Gamma}_p \vecc{\theta}_t dt + \vecc{\Sigma}_p  d\vecc{W}_t, \\
d\vecc{\eta}_t &= -\vecc{\Gamma}_a \vecc{\eta}_t dt + \vecc{\Sigma}_a  d\vecc{U}_t,
\end{align}
where $\vecc{\zeta}_t$ is a mean-zero, stationary Gaussian colored noise with covariance function equals to $\vecc{\kappa}(t) = \vecc{C}_p e^{-\vecc{\Gamma}_p t} \vecc{M}_p \vecc{C}_p^T \in \RR^{d_1 \times d_1}$, and $\vecc{U}_t$, $\vecc{W}_t$ are independent Wiener processes.   

In the absence of $\vecc{\sigma}_a(\vecc{x}_t) \vecc{\eta}_t$, the model can be derived from a microscopic Hamiltonian model describing a particle interacting with an equilibrium heat bath at temperature $T$. Therefore, the above model describes a system driven out of equilibrium by the active force $\vecc{\sigma}_a(\vecc{x}_t) \vecc{\eta}_t$. 
The above model can be viewed as a closely related variant of the ones studied in \citep{leyman2018tuning}. In the joint limit where $\vecc{\kappa}(t)$ tends to a Dirac delta function (memoryless limit), $\vecc{\zeta}_t$ tends to a white noise (white noise limit) and $m \to 0$ (small mass limit), we recover the active Ornstein-Uhlenbeck model for active matter systems studied in \citep{dabelow2019irreversibility} but with inhomogeneous activity due to the state-dependence of $\vecc{\eta}_a$ here.

\ex {\bf A charged particle in a spatially inhomogeneous magnetic field. } \label{eg_3} 
We consider an electrically charged  particle of charge $q$ in an equilibrium homogeneous heat bath. It is subject to a position-dependent magnetic field $\vecc{B}(\vecc{x})$ ($\vecc{x} \in \RR^3$) \citep{vuijk2019anomalous} and time-dependent force field, $\vecc{F} = -\vecc{\nabla}_{\vecc{x}}U(t,\vecc{x}) + q \vecc{E}(t,\vecc{x})$, consisting of forces from conservative potential and electric field.  Assuming that the magnetic field is pointing along the unit vector $\vecc{n}$ and $B(\vecc{x})$ is the magnitude (i.e. $\vecc{B}(\vecc{x}) = B(\vecc{x}) \vecc{n}$),  the Lorentz force $q \vecc{v}_t \times \vecc{B}(\vecc{x}_t)$  can be written as $q B(\vecc{x}_t) \vecc{Z} \vecc{v}_t$, where $\vecc{Z}$ is a matrix with elements given by $Z_{ij} = -\epsilon_{ijk} n_k$, where $\epsilon_{ijk}$ is the totally antisymmetric Levi-Civita symbol in 3D and $n_k$ is the $k$th component of $\vecc{n}$. This system can be described by  the GLE with $\vecc{\gamma}_0(\vecc{x}) = -q B(\vecc{x}) \vecc{Z}$, $\vecc{\sigma}_0 = \vecc{0}$, $\vecc{g} = \vecc{h}^T = \vecc{\sigma}_b$, $\vecc{\sigma} = \sqrt{k_B T} \vecc{\sigma}_b$, $\vecc{\xi}_t$ is the same colored noise as introduced in Section \ref{sect_gles}  but with its covariance function equals to $\vecc{\kappa}(t)$: 
\begin{align}
d\vecc{x}_t &= \vecc{v}_t dt, \\ 
m  d\vecc{v}_t &= \vecc{F}(t,\vecc{x}_t) dt - \vecc{\sigma}_b \left(\int_0^t \vecc{\kappa}(t-s) \vecc{\sigma}_b^T \vecc{v}_s ds \right) dt + q B(\vecc{x}_t) \vecc{Z} \vecc{v}_t  dt + \sqrt{k_B T} \vecc{\sigma}_b \vecc{\xi}_t dt.
\end{align}

In the Markovian limit (i.e. joint memoryless and white noise limit), one obtain a Langevin-Kramers equation with a state-dependent damping term (with a positive stable but not positive definite effective ``damping'' matrix) and an additive white noise term (c.f. \citep{pavliotis2010asymptotic}). The source of the state-dependence in the ``damping'' comes solely from the magnetic field. Different variants of model for such system have  been studied in \citep{hidalgo2016non, lisy2013brownian, Harko2016,gle_oscillatingfields, vuijk2019anomalous,PhysRevE.97.032117}. 

%One could also consider the more general case where $\vecc{\sigma}$ and/or $T$ is state-dependent. 
%small mass limit only is not well-defined. 

\section{Homogenization for a Class of SDEs with State-Dependent Coefficients}
\label{sect_generalhomogthm}

%relax boundness ass up to explosion time?

%homog thoery approach pavlioris/stuart
%other approach?
%backgrounds?
%existence and uniquess of our sdes? hsu

In this section, we recall a homogenization result that will be needed for studying homogenization for our GLEs and their functionals. This result is a special case of the main theorem in \citep{2019Lim}.

%literature review of homogenization 
%papanicol lec?

Let $n_1$, $n_2$, $k_1$, $k_2$ be positive integers. Let $\epsilon > 0$ be a small parameter and $\vecc{X}^{\epsilon}(t) \in \RR^{n_1}$, $\vecc{Y}^{\epsilon}(t) \in \RR^{n_2}$ for $t \in [0,T]$, where  $T>0$ is a constant. Let $\vecc{W}^{(k_1)}$ and $\vecc{W}^{(k_2)}$ denote independent Wiener processes, which are $\RR^{k_1}$-valued and $\RR^{k_2}$-valued respectively, on a filtered probability space $(\Omega, \mathcal{F}, \mathcal{F}_t, \mathbb{P})$ satisfying the usual conditions \citep{karatzas2014brownian}.
%Here $\mathcal{F}_t$ is the usual augmentation of $\sigma(\{\vecc{W}_s := (\vecc{W}^{(k_1)}_s, \vecc{W}^{(k_2)}_s) \in \RR^{k_1+k_2} : s \leq t \})$.

With respect to the standard bases of  $\RR^{n_1}$ and $\RR^{n_2}$ respectively, we write:
\begin{align}
\vecc{X}^{\epsilon}(t) &= (X_1^{\epsilon}(t),X^{\epsilon}_2(t),\dots, X^{\epsilon}_{n_1}(t)),  \\
\vecc{Y}^{\epsilon}(t) &= (Y^{\epsilon}_1(t),Y^{\epsilon}_2(t),\dots, Y^{\epsilon}_{n_2}(t)).
\end{align}

We consider the following family of singularly perturbed SDE systems\footnote{Note that here the variables $\vecc{X}^\epsilon(t)$ and $\vecc{Y}^\epsilon(t)$ are general and they do not necessarily represent position and velocity variables of a physical system.} for 

\noindent $(\vecc{X}^\epsilon(t), \vecc{Y}^\epsilon(t)) \in \RR^{n_1} \times \RR^{n_2}$:
\begin{align}
d\vecc{X}^{\epsilon}(t) &= \vecc{A}_{1}(t,\vecc{X}^{\epsilon}(t)) \vecc{Y}^{\epsilon}(t) dt + \vecc{B}_{1}(t,\vecc{X}^{\epsilon}(t)) dt + \vecc{\Sigma}_{1}(t,\vecc{X}^{\epsilon}(t)) d\vecc{W}^{(k_1)}(t), \label{sde1} \\
\epsilon d\vecc{Y}^{\epsilon}(t) &= \vecc{A}_{2}(t,\vecc{X}^{\epsilon}(t)) \vecc{Y}^{\epsilon}(t) dt + \vecc{B}_{2}(t,\vecc{X}^{\epsilon}(t)) dt + \vecc{\Sigma}_2(t,\vecc{X}^{\epsilon}(t)) d\vecc{W}^{(k_2)}(t), \label{sde2}
\end{align}
with the initial conditions, $\vecc{X}^{\epsilon}(0) = \vecc{X}^\epsilon$ and $\vecc{Y}^{\epsilon}(0) = \vecc{Y}^\epsilon$, where $\vecc{X}^\epsilon$ and $\vecc{Y}^\epsilon$ are random variables that possibly depend on $\epsilon$. In the SDEs \eqref{sde1}-\eqref{sde2}, the coefficients $\vecc{A}_1: \RR^+ \times \RR^{n_1} \to  \RR^{n_1 \times n_2}$, $\vecc{A}_2 : \RR^+ \times \RR^{n_1}  \to   \RR^{n_2 \times n_2}$, $\vecc{\Sigma}_2 : \RR^+ \times \RR^{n_1}  \to \RR^{n_2 \times k_2}$ are non-zero matrix-valued functions, whereas $\vecc{B}_1 : \RR^+ \times  \RR^{n_1} \to \RR^{n_1}$, $\vecc{B}_2 : \RR^+ \times \RR^{n_1} \to \RR^{n_2}$,
$\vecc{\Sigma}_1 : \RR^+ \times \RR^{n_1} \to \RR^{n_1 \times k_1}$ are (possibly zero) matrix-valued or vector-valued functions. They may depend on $\vecc{X}^{\epsilon}$,  as well as on $t$  explicitly, as indicated by the parenthesis $(t, \vecc{X}^{\epsilon}(t))$.  

%We will allow the possibility of explosion of the solutions to \eqref{sde1}-\eqref{sde2} and to \eqref{mainlimitingeqn}, as this will cover a larger class of physically relevant systems (as stronger assumptions on the coefficients like the linear growth condition \citep{karatzas2012brownian} are not made to ensure global existence of solutions). Indeed, in many physical applications the state space of $\vecc{x}^{\epsilon}(t)$ rarely coincides with the whole space $\RR^{n_1}$.
%Therefore, the SDEs considered are general enough to model systems with multiple time scales.

We are interested in the limit as $\epsilon \to 0$ of the SDEs \eqref{sde1}-\eqref{sde2}, in particular the limiting behavior of the  process $\vecc{X}^{\epsilon}(t)$, under  appropriate assumptions\footnote{We forewarn the readers that our assumptions can be relaxed in various directions (see the relevant remarks in \citep{2019Lim}) but we will not pursue these generalizations here. This approach may not be too appealing from a mathematical point of view but we stress that the main goal of the paper is to communicate, in the simplest yet rigorous manner, the consequences of the homogenization results to a broad range of audience and therefore some sacrifices in the completeness are unavoidable.}  on the coefficients. We make the following assumptions concerning the SDEs \eqref{sde1}-\eqref{sde2} and  \eqref{mainlimitingeqn}.

\begin{ass} \label{aexis} The global solutions, defined on $[0,T]$, to the pre-limit SDEs \eqref{sde1}-\eqref{sde2} and to the limiting SDE \eqref{mainlimitingeqn} a.s. exist and are unique for all $\epsilon > 0$ (i.e. there are no explosions).
\end{ass}

\begin{ass} \label{a0_ch2} The matrix-valued functions 
$$\{ -\vecc{A}_2(t,\vecc{X}); t \in [0,T], \vecc{X} \in \RR^{n_1}\}$$ are {\it uniformly positive stable}, i.e. all real parts of the eigenvalues of $-\vecc{A}_2(t, \vecc{X})$ are bounded from below, uniformly in $t$ and $\vecc{X}$, by a positive constant (or, equivalently, the matrix-valued functions $\{\vecc{A}_2(t, \vecc{X});  t \in [0,T], \vecc{X} \in \RR^{n_1}\}$ are {\it uniformly Hurwitz stable}). 
\end{ass}
%do we need jacobian to be continuous?

%This assumption is crucial to obtain desired estimate on the moments of certain processes (see Proposition \ref{expbound}).

%We provide sufficient conditions that guarantee this additional assumption to hold in the appendix.

%Note that this assumption is equivalent to the assumption that the following deterministic Cauchy problem is uniformly well-posed: \begin{equation} \vecc{u}'(t) = \vecc{A}_2(\vecc{x}_t,\epsilon) \vecc{u}(t), \ \ \vecc{u}(0) = \vecc{u} \in \mathcal{H}^{\vecc{v}}. \end{equation}
%In particular, this assumption implies that there exists constants $M > 0$ and $\alpha \geq 0$ such that $$\|\vecc{S}(t)\| \leq M e^{\alpha t}$$ for $t \geq 0$.

%However, we need an decay estimate on $\|\vecc{S}(t)\|$ to prove the main result.

%$\alpha$ is the (exponential) growth bound of the semigroup, defined by $$\alpha = \alpha(\vecc{U}) := \inf\{ w \in \RR : \text{ there exists } M_{w} \geq 1 \text{ such that }  \|\vecc{S}(t)\| \leq M_{w} e^{wt} \text{ for } t \in [0,T]\}.$$ As we will see, the negativity of the growth bound, $\alpha$, is crucial to obtain desired estimate on the moments of certain processes. We provide sufficient conditions that guarantee this additional assumption to hold in the appendix.

%In particular, this assumption implies that the Cauchy problem is exponentially stable, i.e. its solution, $\vecc{S}(\cdot)$, admits an exponential bound: there exists constants $\alpha \geq 0$ and $M \geq 1$ sich that

\begin{ass} \label{a1_ch2} For $t \in [0,T]$, $\vecc{X} \in \RR^{n_1}$, and $i=1,2$, the functions $\vecc{B}_i(t,\vecc{X})$ and
$\vecc{\Sigma}_i(t,\vecc{X})$ are continuous and bounded in $t$ and $\vecc{X}$, and Lipschitz in $\vecc{X}$, whereas the functions $\vecc{A}_i(t,\vecc{X})$ and $(\vecc{A}_i)_{\vecc{X}}(t,\vecc{X})$ are continuous in $t$, continuously differentiable in $\vecc{X}$, bounded in $t$ and $\vecc{X}$, and Lipschitz in $\vecc{X}$.
Moreover, the  functions $(\vecc{A}_i)_{\vecc{X} \vecc{X}}(t,\vecc{X})$ ($i=1,2$) are bounded for every $t \in [0,T]$ and $\vecc{X} \in \RR^{n_1}$.
\end{ass}

\begin{ass} \label{a2_ch2} The initial condition $\vecc{X}^\epsilon(0) = \vecc{X}^\epsilon \in \RR^{n_1}$ is an $\mathcal{F}_0$-measurable random variable that may depend on $\epsilon$, and we assume that $\mathbb{E}[|\vecc{X}^\epsilon|^p] = O(1)$ as $\epsilon \to 0$ for all $p>0$. Also, $\vecc{X}^\epsilon$ converges, in the limit as $\epsilon \to 0$, to a random variable $\vecc{X}$ as follows:
$\mathbb{E}\left[|\vecc{X}^\epsilon - \vecc{X}|^p \right]  = O(\epsilon^{p r_0})$, where $r_0 > 1/2$ is a constant, as $\epsilon \to 0$.       The initial condition $\vecc{Y}^\epsilon(0) = \vecc{Y}^\epsilon \in \RR^{n_2}$ is an $\mathcal{F}_0$-measurable random variable that may depend on $\epsilon$, and we assume that for every $p>0$, $\mathbb{E}[ |\epsilon \vecc{Y}^\epsilon|^p] = O(\epsilon^\alpha)$ as $\epsilon \to 0$, for some $\alpha \geq p/2$.

%Both $\vecc{x}$ and $\vecc{v}^\epsilon$ (for every $\epsilon > 0$) are independent of the $\sigma$-algebra generated by $\vecc{W}$.

%The initial conditions, $\vecc{x} \in \RR^{n_1}$, $\vecc{v} \in \RR^{n_2}$, are $\mathcal{F}_0$-measurable random variables independent of $\epsilon$ and have finite moments of all orders, i.e.
%$\mathbb{E} |\vecc{x}|^p,  \ \mathbb{E} |\vecc{v}|^p < \infty$ for all $p > 0$.

%for every $t \in [0,T]$, $\epsilon >0$, the solutions to the pre-limit equations \eqref{sde1}-\eqref{sde2} have finite moments of all orders up to the explosion time, i.e. for $t \in [0,\tau)$, $\mathbb{E} |\vecc{x}^{\epsilon}(t)|^p,  \ \mathbb{E} |\vecc{v}^{\epsilon}(t)|^{p} < \infty$ for all $p > 0$.
\end{ass}

We now state the homogenization theorem.

\begin{thm} \label{mainthm}
Suppose that the family of SDE systems $\eqref{sde1}$-$\eqref{sde2}$ satisfies Assumption \ref{aexis}-\ref{a2_ch2}. Let $(\vecc{X}^{\epsilon}(t), \vecc{Y}^{\epsilon}(t)) \in \RR^{n_1} \times \RR^{n_2}$ be their solutions, with the initial conditions $(\vecc{X}^\epsilon, \vecc{Y}^\epsilon)$. Let $\vecc{X}(t) \in \RR^{n_1}$ be the solution to the following It\^o SDE with the initial position $\vecc{X}(0) = \vecc{X}^\epsilon$:
\begin{align}
d\vecc{X}(t) &= [\vecc{B}_1(t,\vecc{X}(t))-\vecc{A}_1(t,\vecc{X}(t))\vecc{A}_2^{-1}(t,\vecc{X}(t))\vecc{B}_2(t,\vecc{X}(t))] dt \nonumber \\
&\ \ \ \  + \vecc{S}(t,\vecc{X}(t)) dt  + \vecc{\Sigma}_1(t,\vecc{X}(t)) d\vecc{W}^{(k_1)}(t) \nonumber \\
&\ \ \ \  - \vecc{A}_1(t,\vecc{X}(t)) \vecc{A}_2^{-1}(t,\vecc{X}(t))\vecc{\Sigma}_2(t,\vecc{X}(t)) d\vecc{W}^{(k_2)}(t). \label{mainlimitingeqn}
\end{align}
In the above  $\vecc{S}(t,\vecc{X}(t))$ is the {\it noise-induced drift vector} whose $i$th component is given by
\begin{equation} \label{non_indexfree}
S^{i}(t,\vecc{X}) = -\frac{\partial}{\partial X^{l}} \bigg((A_1 A_2^{-1})^{ij}(t,\vecc{X}) \bigg) \cdot A_1^{lk}(t,\vecc{X}) \cdot J^{jk}(t,\vecc{X}), 
\end{equation}
where $i,l=1,\dots,n_1, \  j,k=1,\dots,n_2$, or in index-free notation, 
\begin{equation} \label{indexfree}
\vecc{S} = \vecc{A}_1 \vecc{A}_2^{-1} \vecc{\nabla}\cdot (\vecc{J}\vecc{A}_1^T) -\vecc{\nabla} \cdot (\vecc{A}_1 \vecc{A}_2^{-1} \vecc{J} \vecc{A}_1^T) , 
\end{equation}
and $\vecc{J} \in \RR^{n_2 \times n_2}$ is the unique solution to the Lyapunov equation:
\begin{equation} \label{lyp}
\vecc{J} \vecc{A}_2^{T} + \vecc{A}_2 \vecc{J} = -\vecc{\Sigma}_2  \vecc{\Sigma}_2^{T}.
\end{equation}
Then the process $\vecc{X}^{\epsilon}(t)$ converges, as $\epsilon \to 0$, to the solution $\vecc{X}(t)$, of the It\^o SDE \eqref{mainlimitingeqn}, in the following sense: for all finite $T>0$,
\begin{equation}
 \sup_{t \in [0,T]} |\vecc{X}^\epsilon(t) - \vecc{X}(t)| \to 0, \end{equation}
in probability, in the limit as $\epsilon \to 0$.
\end{thm}

\begin{rmk}
If $\vecc{\Sigma}_1$ and $\vecc{\Sigma}_2$ are independent of $\vecc{X}$, then the It\^o equation \eqref{mainlimitingeqn} is equivalent to the  equation:
\begin{align}
d\vecc{X}(t) &= [\vecc{B}_1(t,\vecc{X}(t))-\vecc{A}_1(t,\vecc{X}(t))\vecc{A}_2^{-1}(t,\vecc{X}(t))\vecc{B}_2(t,\vecc{X}(t))] dt + \vecc{H}_\alpha(t,\vecc{X}(t)) dt  \nonumber \\
&\ \ \ \   + \vecc{\Sigma}_1(t) d\vecc{W}^{(k_1)}(t)  - \vecc{A}_1(t,\vecc{X}(t)) \vecc{A}_2^{-1}(t,\vecc{X}(t))\vecc{\Sigma}_2(t) \circ^{\alpha} d\vecc{W}^{(k_2)}(t), \label{mainlimitingeqn_strat}
\end{align}
where $\circ^{\alpha}$, $\alpha \in [0,1]$, specifies the rule of stochastic integration, whereby the stochastic integral is evaluated at $t_n = (1-\alpha)t_n + \alpha t_{n+1}$ on the discretization intervals $[t_n, t_{n+1}]$ (so $\alpha=0$ corresponds to It\^o integral, $\alpha=1/2$ to Stratonovich, and $\alpha=1$ to anti-It\^o), and $\vecc{H}_{\alpha}$ is the corresponding noise-induced drift term whose $i$th component is:
\begin{align}
H_{\alpha}^i &= S^i - \alpha \frac{\partial (A_1 A_2^{-1} \Sigma_2)^{ik}}{\partial X^j} (A_1 A_2^{-1} \Sigma_2)^{jk},
\end{align}
with $S^i$ given by \eqref{non_indexfree}. 

After some algebraic manipulations and using the Lyapunov equation $\vecc{A}_2 \vecc{J} + \vecc{J} \vecc{A}_2^T = - \vecc{\Sigma}_2 \vecc{\Sigma}_2^T$, one can rewrite $H^i_\alpha$ as: 
\begin{equation}
H^i_\alpha = \frac{1}{2}Q^{qj}(\alpha) [G_q, G_j]^i, 
\end{equation}
where $\vecc{G}_q$ denotes the vector field associated to the $q$th column of the matrix $\vecc{A}_1 \vecc{A}_2^{-1}$,  $[G_q, G_j]^i$ denotes the $i$th component of Lie bracket\footnote{If $\vecc{A}$ and $\vecc{B}$ are the first order differential operators corresponding to the vector fields $\vecc{A}(\vecc{x})$ and $\vecc{B}(\vecc{x})$, i.e. $\vecc{A}= \sum_{i} A^i(\vecc{x}) \frac{\partial}{\partial x^i}$ and $\vecc{B} = \sum_j B^j(\vecc{x}) \frac{\partial}{\partial x^j}$, then the Lie bracket (commutator) between $\vecc{A}$ and $\vecc{B}$ is defined as the operator $[\vecc{A}, \vecc{B}] = \vecc{A} \vecc{B} - \vecc{B}\vecc{A}$.} of the vector fields $\vecc{G}_q$ and $\vecc{G}_j$ (i.e. the derivative of $\vecc{G}_j$ along the flow generated by $\vecc{G}_q$), and
\begin{equation}
\vecc{Q}(\alpha) = \alpha \vecc{J}\vecc{A}_2^T - (1-\alpha) \vecc{A}_2 \vecc{J}.
\end{equation}
Provided that $\vecc{A}_2$ is Hurwitz stable, $\vecc{Q}(\alpha)$ can be represented as the solution to the Lyapunov equation \citep{bellman1997introduction}:
\begin{equation}
\vecc{A}_2 \vecc{Q}(\alpha) + \vecc{Q}(\alpha) \vecc{A}_2^T = antisym(\vecc{A}_2 \vecc{\Sigma}_2 \vecc{\Sigma}_2^T),
\end{equation}
where $antisym(\vecc{A})$ denotes the antisymmetric part of the matrix $\vecc{A}$.

Now, let us consider the Stratonovich case $\alpha=1/2$. In this case, $\vecc{Q} := \vecc{Q}(1/2)$ is the antisymmetric part of the Onsager matrix $-\vecc{A}_2 \vecc{J}$, i.e. $\vecc{Q} = (\vecc{J} \vecc{A}_2^T - \vecc{A}_2 \vecc{J})/2$ (see also Remark \ref{rmk41}). Therefore, when the {\it detailed balance condition} (i.e. when $\vecc{A}_2 \vecc{\Sigma}_2 \vecc{\Sigma}_2^T$ is symmetric) holds, $\vecc{Q}$ (physically a measure of irreversibility of the fast process, and mathematically a matrix encoding stochastic area of the limiting process) vanishes and the resulting limiting SDE for $\vecc{X}(t)$ is a Stratonovich SDE without additional drift correction terms. On the other hand, if $\alpha=0$ (It\^o), $\vecc{Q}(\alpha=0)$ is simply the (non-zero) Onsager matrix, whereas if $\alpha=1$ (anti-It\^o), $\vecc{Q}(\alpha=1)$ equals to negative transpose of the Onsager matrix.  
\end{rmk}

%renormalized SDEs, rough paths? 
%connection to theory of rough path theory and regularity structures

\end{document}